\documentclass[]{interact}
\usepackage[font=small,labelsep=quad,labelfont=bf]{caption}
\usepackage{url}
\usepackage{algorithm, algorithmic}
\usepackage{multirow}
\usepackage{mathrsfs,amsmath,amsfonts,amssymb}
\usepackage{color}
\usepackage{subfigure}
\usepackage{MnSymbol}
\usepackage{threeparttable}
\usepackage{hyperref}
\hypersetup{
    colorlinks=true,
    citecolor=black,
    linkcolor=black,
    filecolor=black,
    urlcolor=black,
}
\usepackage{epstopdf}

\usepackage[numbers,sort&compress]{natbib}
\theoremstyle{plain}
\newtheorem{theorem}{Theorem}[section]
\newtheorem{proposition}[theorem]{Proposition}
\newtheorem{condition}{Condition}

\theoremstyle{definition}

\theoremstyle{remark}

\def\hat{\widehat}

\def\bar{\overline}
\def\epsilon{\varepsilon}

\def\dfrac{\displaystyle\frac}

\def\*{$\!\!^{^{^{\displaystyle *}}}$}

\def\Rmn#1{\expandafter\uppercase\expandafter{\romannumeral #1}}

\begin{document}

\title{Functional sufficient dimension reduction through distance covariance}
\author{
\name{Xing Yang\textsuperscript{a}, Jianjun Xu\textsuperscript{b*}
\thanks{*Corresponding authors: Jianjun Xu. E-mail: xjj1994@ustc.edu.cn}}
\affil{ \textsuperscript{a}Department of Mathematics, Northeastern University, Boston, MA 02115, USA;
\textsuperscript{b}International Institute of Finance, School of Management, University of Science and Technology
of China, Hefei, 230026, Anhui, China.}
}

\maketitle

\begin{abstract}
Our research proposes a novel method for reducing the dimensionality of functional data, specifically for the case where the response is a scalar and the predictor is a random function. Our method utilizes distance covariance, and has several advantages over existing methods. Unlike current techniques which require restrictive assumptions such as linear conditional mean and constant covariance, our method has mild requirements on the predictor. Additionally, our method does not involve the use of the unbounded inverse of the covariance operator.
The link function between the response and predictor can be arbitrary, and our proposed method maintains the advantage of being model-free, without the need to estimate the link function. Furthermore, our method is naturally suited for sparse longitudinal data.
We utilize functional principal component analysis with truncation as a regularization mechanism in the development of our method. We provide justification for the validity of our proposed method, and establish statistical consistency of the estimator under certain regularization conditions.
To demonstrate the effectiveness of our proposed method, we conduct simulation studies and real data analysis. The results show improved performance compared to existing methods.
\end{abstract}

\section{Introduction}
In contemporary data analysis, functional data are  prevalent in many applications such as speech recognition, magnetic resonance imaging (MRI), online handwriting recognition and longitudinal data analysis \citep{ramsay1991some}.
Under a functional data analysis (FDA) framework, each sample element is considered to be a function.
A hot issue is to study how a response variable varies with a random function $X(t)$, where $t$ is an index variable defined on an interval.

\begin{figure}[t]
  \centering
  \includegraphics[width=7cm]{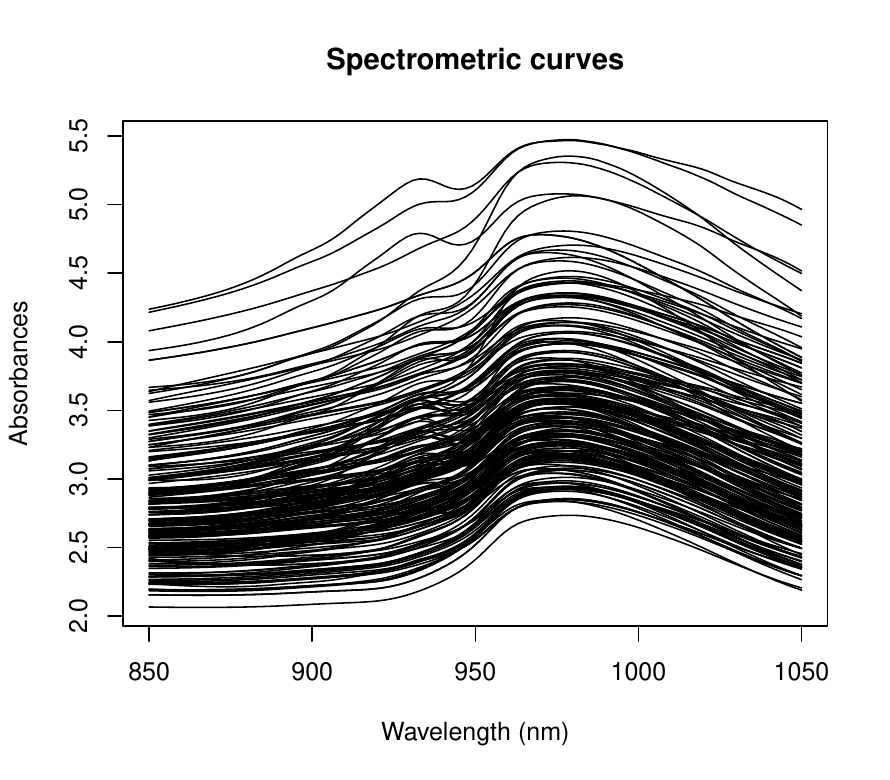}
  \caption{The near infrared
absorbance spectrum curves.}\label{spectrum}
\end{figure}

Take Tecator data as an example to introduce the problem of functional regression.
Each sample of this dataset contains finely chopped pure meat with different moisture,
fat and protein contents. Using analytical chemistry to measure fat content is expensive,
while infrared analysis is substantially cheaper.
The aim of the analysis is predicting the fat content of pieces of meat from a near infrared
absorbance spectrum which is a curve, see Figure \ref{spectrum}.
This is a typical functional regression problem which has been investigated from both parametric and nonparametric point of views \citep{cardot1999functional,ferraty2002functional}. However, parametric modeling can be restrictive  in some applications, while nonparametric modeling can be unworkable due to the infinite-dimensionality of the functional data.
For instance, most functional regression and correlation measure problems involve the inverse of compact operators which are unbounded.
This is triggered by the infinite-dimensionality of the functional data, therefore,  dimension reduction is key for functional data
modeling and analysis.

Despite the infinite-dimensional nature of functional data, interestingly, the data set tends to have a certain pattern which might be represented by finite indexes.
With the rapid development of
functional data analysis,  functional sufficient dimension reduction (FSDR) problems have received increasing attention in the literature. \cite{ferre2003functional} first extended sliced inverse regression (SIR, \citealp{li1991sliced}) to the functional case and assumed
\begin{equation}\label{fSDR}
  Y=g\left(\left\langle\beta_{1}, X\right\rangle, \ldots,\left\langle\beta_{K}, X\right\rangle, \epsilon\right),
\end{equation}
where the response $Y$ is a random variable, the predictor $X$ takes values in a functional space $\mathcal{H}$ with inner product denoted by
$\left\langle\cdot,\cdot \right\rangle$,~ $\beta_1,\ldots, \beta_K$ are $K$ linearly independent functions in $\mathcal{H}$, $\epsilon$ is a scalar random noise, and $g$ is an unknown function from $\mathbb{R}^{K+1}$ to $\mathbb{R}$.
In the FSDR literature, the subspace spanned by $\beta_1,\ldots,\beta_K$ is called the functional sufficient dimension reduction  subspace. After \cite{ferre2003functional} proposed the functional SIR (FSIR),
quite a few different methods have been developed for estimating the SDR space. For instance,
functional inverse regression \citep{ferre2005smoothed}, functional contour regression \citep{wang2013functional},
functional K-means inverse regression \citep{wang2014functional}, functional sliced average variance estimation
(FSAVE, \citealp{lian2014series}), the hybrid method of FSIR and FSAVE \citep{wang2015hybrid}, functional cumulative slicing \citep{yao2015effective}, robust FSIR \citep{wang2017robust} and functional directional regression \citep{wang2020directional}.
\cite{lian2015functional} consider FSIR and FSAVE via a Tikhonov regularization approach and show that their convergence rates are the same as the
minimax rates for functional linear regression.
\cite{li2017nonlinear} proposed functional generalized SIR and functional generalized SAVE for nonlinear
sufficient dimension reduction where both the predictor and the response may
be random functions.
\cite{li2022dimension} and \cite{lee2022functional} developed sufficient dimension reduction methods for function-on-function
regression through weak conditional moments and average Fr{\'e}chet derivatives respectively.
\cite{hsing2009rkhs} proposed a method under reproducing kernel Hilbert space (RKHS) framework which can be applied to
finite or infinite dimensional predictor space in a unified framework.
Other researches include but are not limited to localized and regularized
versions of FSIR and FSAVE \citep{wang2018functional} and functional surrogate assisted slicing \citep{wang2019dimension}
which are aimed at SDR for binary classification. More information can be referred to the review article \cite{song2019sufficient} and the monograph \cite{li2018sufficient}.
However all of the above methods need the linear conditional mean assumption or constant covariance assumption or
both and these assumptions are not easy to verify in practice. Most of the above methods need to estimate the inverse of a covariance operator which is unbounded since the covariance operator is defined on an infinite-dimensional space.

Existing literature on functional sufficient  dimension reduction is mainly based on the sliced inverse moment methods.
However, take FSIR \citep{ferre2003functional}
 as an example, applying FSIR to sparsely observed longitudinal data is practically infeasible, since  choosing a sufficiently large number
of slices would result in too few observations in each slice with which to estimate a conditional covariance operator.
Therefore, there is very little literature on sufficient dimension reduction for sparse longitudinal data. \cite{jiang2014inverse} extended the method of \cite{ferre2005smoothed} for sparse longitudinal data.
Functional cumulative slicing \citep{yao2015effective} and some other
methods are also suitable for longitudinal data. We will show that the proposed  method is applicable to both dense functional data and sparse longitudinal data.

In multivariate setting, \cite{sheng2016sufficient} proposed a method for SDR via distance covariance \citep{szekely2007measuring,szekely2009brownian,lyons2013distance}. Inspired by this method, we extend it to the functional context which has not been considered before.
The goal of this paper is to develop a class of sufficient dimension reduction techniques for functional data that require
no inversion of the covariance operator, using the idea of distance covariance. To the best of our knowledge, this is the first time
that distance covariance methodology is extended beyond the usual multivariate regression setting to functional data analysis. An
important contribution of this paper is to bridge the gap between the nascent area of dependence measure, functional data
analysis, and sufficient dimension reduction.

In this article, following the work of \cite{sheng2016sufficient}, we first
use distance covariance for functional sufficient dimension reduction.
This method does not require  linear conditional mean assumption and constant covariance assumption. It also does not involve the inverse of the covariance operator which is not bounded.
Under mild conditions, we prove the
validity of the proposed  method as a sufficient functional dimension reduction method and we use functional principal components
method as a form of regularization to make it feasible to estimate a infinite-dimension function in a finite-dimensional subspace.
We also construct the consistency of the proposed estimator.
Simulation and real data analysis are conducted to exhibit the superiority of the proposed method.

The rest of the paper is organized as follows. In Section
\ref{methodology}, we introduce distance covariance, functional sufficient dimension reduction and propose our method for
functional sufficient dimension reduction via distance covariance at the population level. Finite-sample estimation and its statistical consistency are presented in Section \ref{estmation}.
Simulations and  real data analysis are carried out in Section \ref{numerical}. Section \ref{condis} concludes the paper, and all the proofs are deferred to the Appendix.

\section{Methodology}\label{methodology}
\subsection{Distance covariance}\label{DC}
Distance covariance (DCOV) proposed by \cite{szekely2007measuring} is a new measure of dependence between random vectors. The appealing property of distance covariance is that it is zero if and only if the random variables are independent.
In this subsection $U$ in $\mathbb{R}^p$ and $V$ in $\mathbb{R}^q$ are random vectors, where $p$
and $q$ are positive integers. The Euclidean norm of $x$ in $\mathbb{R}^p$ is $|x|_p$. The characteristic functions of $U$, $V$ and $(U, V)$ are denoted by $f_U$, $f_V$ and $f_{UV}$ respectively. Then the DCOV defined in \cite{szekely2007measuring} is the nonnegative number
$\mathcal{V}(U,V)$ with
\begin{equation}\label{DCOV}
  \mathcal{V}^{2}(U, V)=\int_{\mathbb{R}^{2}}\left|f_{U V}(s, t)-f_{U}(s) f_{V}(t)\right|^{2} w(s, t) d s d t,
\end{equation}
where $|f|^2=f\cdot \bar{f}$ and $w(s,t)$ is a weight function. If we choose $w(s,t)= (\pi^2s^2t^2)^{-1}$, \cite{szekely2009brownian} gave an equivalent form of DCOV as
\begin{equation}\begin{aligned}
\mathcal{V}^{2}(U, V) &=E|U-U^{\prime} |_p| V-V^{\prime}|_q+E|U-U^{\prime}|_p E|V-V^{\prime}|_q \\
&-E|U-U^{\prime}|_p|V-V^{\prime \prime}|_q-E| U-U^{\prime \prime}|_p |V-V^{\prime}|_q,
\end{aligned}\end{equation}
where $(U,V)$, $(U^{\prime},V^{\prime})$ and $(U^{\prime\prime},V^{\prime\prime})$ are i.i.d.
Here we list several useful properties of DCOV. For random vectors $U\in \mathbb{R}^p$ and $V\in \mathbb{R}^q$ such that $E(|U|_p+|V|_q)<\infty$, the following properties hold:
\begin{enumerate}
  \item[(i)] $\mathcal{V}(U, V)=0$ if and only if $U$ and $V$ are independent.
  \item[(ii)] $\mathcal{V}(a_1+b_1C_1U,a_2+b_2C_2V)=\sqrt{|b_1b_2|}\mathcal{V}(U,V)$ for all constant vectors $a_1\in \mathbb{R}^p$, $a_2\in \mathbb{R}^q$, scalars $b_1$, $b_2$ and orthonormal matrix $C_1$, $C_2$ in $\mathbb{R}^{p\times p}$ and $\mathbb{R}^{q\times q}$, respectively.
  \item[(iii)] If the random vector $(U_1, V_1)$ is independent of the random vector $(U_2, V_2)$, then
  \begin{equation*}
    \mathcal{V}(U_1+U_2,V_1+V_2)\leq \mathcal{V}(U_1,V_1)+\mathcal{V}(U_2,V_2).
  \end{equation*}
  Equality holds if and only if $U_1$ and $V_1$ are both constants, or $U_2$ and $V_2$ are both
constants, or $U_1$, $U_2$, $V_1$, $V_2$ are mutually independent.
\end{enumerate}

As mentioned in \cite{sheng2016sufficient} and \cite{zhang2019robust}, property (i) makes it possible that
DCOV can be used as a sufficient dimension reduction tool. The properties (ii) and (iii) will be applied in the subsequent text.
\subsection{Functional sufficient dimension reduction}
The functional sufficient dimension reduction (FSDR) is characterized by conditional independence
\begin{equation}\label{FSDR}
  Y \perp \!\!\!\!\! \perp X \mid\left(\left\langle \beta_{1}, X\right\rangle, \ldots,\left\langle \beta_{K}, X\right\rangle\right),
\end{equation}
where the response $Y$ is a scalar random variable, the predictor $X$ takes values in a functional space $\mathcal{H}$,
$\left\langle\cdot,\cdot \right\rangle$ represents the inner product in $\mathcal{H}$,\ \
$\beta_1,\ldots, \beta_K$ are $K$ linearly independent vectors in $\mathcal{H}$ and $\perp \!\!\!\!\! \perp$ indicates independence.
The multi-index model \eqref{fSDR} is included in \eqref{FSDR}.
Without loss of generality, we consider
$\mathcal{H}=L_2([0,1])$, the space spanned by all the square integrable functions on $[0,1]$. For any $f, g\in L_2([0,1])$,
the inner product is defined by $\left\langle f,g\right\rangle=\int_{0}^{1}f(t)g(t)dt$. The FSDR subspace is denoted by $S=\mathrm{span}\{\beta_1,\ldots,\beta_K\}$.  Obviously, FSDR subspace is not unique, so we only consider the smallest FSDR subspace, which can be defined as the intersection of all FSDR subspace. Following the convention of \cite{wang2019dimension}, we call it the functional central subspace, denoted as $S_{y|x}$. Throughout the article, we assume $S_{y|x}$
exists, which is unique. Let $K$ denote the dimension of $S_{y|x}$. Our primary goal is to   identify $S_{y|x}$ by estimating
$K$ basis functions that span $S_{y|x}$.

\subsection{DCOV as a FSDR method}
We assume $E(X)=\mu_X$ and denote the covariance operator of $X$ by $\Sigma_X=E[(X-\mu_X)\otimes(X-\mu_X)]$, where $\otimes$ is defined as $(f\otimes g)v=\langle f,v\rangle g$ for any $f,g,v\in L_2([0,1])$. For the sequel, we define another inner product on $L_2([0,1])$. For any positive definite self-adjoint and linear operator A, the inner product
$\langle f,g \rangle_A$, for any $f,g \in L_2([0,1])$, is defined as  $\langle f,g \rangle_A=\langle Af, g\rangle=
\langle f, Ag\rangle=\int\int f(s)A(s,t)g(t)dsdt$.
A vector $\bm{f}=(f_1,\ldots,f_{d})$ has $d$ components, where each $f_i$ is a function in $L_2([0,1])$.
We define $\langle \bm{f},\bm{f}\rangle_A=\big(\langle f_i,f_j\rangle_{A}\big)_{1\leq i,j\leq d}$ as a $d\times d$ matrix.
For a function $g\in L_2([0,1])$,
we define $\langle \bm{f},g \rangle=(\langle f_1,g \rangle,\ldots,\langle f_d,g \rangle)^T$ as a $d$-dimensional vector.

Denote $\bm{\beta}=(\beta_1,\ldots,\beta_K)$ where each $\beta_i$ is a function in $L_2([0,1])$.
We will show that under mild conditions, a basis of $S_{y|x}$ can be obtained by solving the following optimization problem:
\begin{equation}\label{maxi}
  \max_{\langle\bm{\beta} , \bm{\beta}\rangle_{\Sigma_X}=I_{K}} \mathcal{V}^2(\langle \bm{\beta},X\rangle,Y),
\end{equation}
where $\langle\bm{\beta} , \bm{\beta}\rangle_{\Sigma_X}=\big(\langle \beta_i,\beta_j\rangle_{\Sigma_X}\big)_{1\leq i,j\leq K}$ is a $K\times K$ matrix, $I_{K}$ is the $K\times K$ identity matrix and
$\langle \bm{\beta},X\rangle=(\langle \beta_1,X\rangle,\ldots,\langle \beta_{K},X\rangle)^T$ is a column vector.
Here we need a scale constraint $\langle\bm{\beta} , \bm{\beta}\rangle_{\Sigma_X}=I_{K}$ to make the maximization procedure work. The reason is that  $\mathcal{V}^2(\langle c\bm{\beta},X\rangle,Y)=|c| \mathcal{V}^2(\langle \bm{\beta},X\rangle,Y)$
for any constant $c$, so we can always get a bigger value of $\mathcal{V}^2(\langle \bm{\beta},X\rangle,Y)$ by multiplying $\bm{\beta}$ by a constant with bigger absolute value.

The following propositions guarantees the validity of DCOV as a tool of functional sufficient dimension reduction.
The solution of the optimization problem \eqref{maxi} indeed spans the functional central subspace.
\begin{proposition}\label{prop1}
Let $\bm{\eta}=(\eta_1, \ldots,\eta_K)$ be a basis of $S_{y|x}$ with $\langle\bm{\eta}, \bm{\eta}\rangle_{\Sigma_X}=I_K$,
$\bm{\beta}=(\beta_1,\ldots, \beta_{K_1})$ with $K_1\leq K$ and $\langle\bm{\beta} , \bm{\beta}\rangle_{\Sigma_X}=I_{K_1}$. Assume $\mathrm{span}(\bm{\beta})\subseteq \mathrm{span}(\bm{\eta})$, then $\mathcal{V}^2(\langle \bm{\beta},X\rangle,Y)\leq \mathcal{V}^2(\langle \bm{\eta},X\rangle,Y)$. The equality holds if and only if $\mathrm{span}(\bm{\beta})=\mathrm{span}(\bm{\eta})$.
\end{proposition}

Proposition \ref{prop1} means that the DCOV between $\langle\bm{\beta}, X\rangle$ and $Y$ is always no more than the DCOV between $\langle\bm{\eta}, X\rangle$ and $Y$ when $\mathrm{span}(\bm{\beta})$ is a subspace of the functional central subspace $\mathrm{span}(\bm{\eta})=S_{y|x}$.
The equality holds if and only if $\mathrm{span}(\bm{\beta})=\mathrm{span}(\bm{\eta})$. However, this result is not enough to guarantee
 DCOV as a tool of functional sufficient dimension reduction. We also  need to consider the situation that  $\mathrm{span}(\bm{\beta})\nsubseteq
 \mathrm{span}(\bm{\eta})$. The next proposition gives the result of this situation under a mild condition.

\begin{condition}\label{assump1}
 Let $\bm{\eta}=(\eta_1, \ldots,\eta_K)$ be a basis of the $S_{y|x}$.
 Denote $\mathrm{span}(\bm{\eta})^{\perp}$ the orthogonal complement space of $\mathrm{span}(\bm{\eta})$ with respect to the inner product $\langle\cdot ,\cdot \rangle_{\Sigma_X}$. We assume that for any $\bm{f}=(f_1,\ldots, f_I)$, $\bm{g}=(g_1,\ldots, g_J)$ where $f_i\in \mathrm{span}(\bm{\eta})$ and $g_j\in \mathrm{span}(\bm{\eta})^{\perp}$, $i=1,\ldots, I$, $j=1,\ldots, J$, we have
 $\langle \bm{f},X\rangle \perp
  \langle \bm{g},X\rangle$.
\end{condition}

This condition is not as strong as it seems to be. When $X$ is a Gaussian process, the independence condition $\langle \bm{f},X\rangle \perp \langle \bm{g},X\rangle$ will be satisfied, because
$Cov(\langle \bm{f},X\rangle, \langle \bm{g},X\rangle)=\bm{0}$. However, Gaussianity is not necessary. Condition \ref{assump1} asymptotically holds when the dimension of $X$ gets reasonably high, see  \cite{sheng2013direction,sheng2016sufficient} for details.

\begin{proposition}\label{prop2}
Let $\bm{\eta}=(\eta_1, \ldots,\eta_K)$ be a basis of $S_{y|x}$ with $\langle\bm{\eta} , \bm{\eta}\rangle_{\Sigma_X}=I_K$,
$\bm{\beta}=(\beta_1,\ldots, \beta_{K_2})$ with $\langle\bm{\beta} , \bm{\beta}\rangle_{\Sigma_X}=I_{K_2}$.  Here $K_2$ could be bigger, less, or equal to $K$.
Assume Condition \ref{assump1} holds and  $\mathrm{span}(\bm{\beta})\nsubseteq \mathrm{span}(\bm{\eta})$, then $\mathcal{V}^2(\langle \bm{\beta},X\rangle,Y)< \mathcal{V}^2(\langle \bm{\eta},X\rangle,Y)$.
\end{proposition}

Proposition \ref{prop2} indicates that if  $\mathrm{span}(\bm{\beta})\nsubseteq \mathrm{span}(\bm{\eta})$, the DCOV between $\langle\bm{\beta}, X\rangle$ and $Y$ is always less than the DCOV between $\langle\bm{\eta}, X\rangle$ and $Y$.
Propositions \ref{prop1} and \ref{prop2} in this article are extensions of propositions 1 and 2 in \cite{sheng2016sufficient}.
Propositions \ref{prop1} and \ref{prop2} guarantee the validity of DCOV as a functional sufficient dimension reduction method,
so we can obtain a basis of $S_{y|x}$ by solving the optimization problem  \eqref{maxi}. Note that the optimization problem  \eqref{maxi} is on the population level and the infinite dimensional functions $X$ and $\bm{\beta}$ make the problem more complicated. Take into account of these problems, some form of regularization is needed in the estimation procedure.  We will elaborate
it in detail in the next section.

\section{Estimation}\label{estmation}
In the previous section, we have establish the method of functional SDR via distance covariance at the population level.
In this section, we will give algorithms for both completely observed functional data and sparse longitudinal data
at the sample level. The
structural dimension $K$ is assumed to be known in this section.

\subsection{Estimation of $S_{y|x}$ for functional data}

In functional context, the estimation procedure is more intractable than that in multivariate context because of the infinite dimensionality of the predictor $X$ and the parameters $\beta_k$, $k=1,\ldots,K$. Practically feasible approaches must include some form of dimensionality reduction.
A standard method in functional data analysis is to embed the infinite dimensional curves into a finite dimensional space. Specifically, $X$ and $\beta_k$ are approximated using the series expansion method. In this article, we consider functional principal components (FPC) basis, which is a common choice in practice \citep{ramsay1991some,yao2005functional,hall2006properties,wang2016functional}.

In reality, the predictor trajectories are observed  intermittently.
For densely observed $X$, individual smoothing can be
used as a pre-processing step to recover smooth trajectories, and the error introduced
by individual smoothing has be shown to be asymptotically negligible under certain design conditions \citep{hall2006properties}. Thus,
to simplify the notation, we only consider completely observed functional data.

Given i.i.d sample $(\bm{X},\bm{Y})=\{(X_i(t),Y_i),\ i=1,2,\ldots,n\}$ and let $\bar{X}(t)=1/n\sum_{i=1}^{n}X_i(t)$, the sample covariance of $X$ can be estimated by $\widehat{\Sigma}_X(s,t)=1/n\sum_{i=1}^{n}\big(X_i(s)-\bar{X}(s)\big)\big(X_i(t)-\bar{X}(t)\big)$.
By Mercer's theorem \citep{riesz1955b}, $\Sigma_X$ and $\widehat{\Sigma}_X$ admit the following eigen-decomposition:
\begin{equation}\label{mercer}
  \Sigma_X(s,t)=\sum_{i=1}^{\infty}\lambda_i\phi_i(s)\phi_i(t),\qquad \widehat{\Sigma}_X(s,t)=\sum_{i=1}^{\infty}\widehat{\lambda}_i\widehat{\phi}_i(s)\widehat{\phi}_i(t),
\end{equation}
where $\{\lambda_j : j\geq1\}$ and $\{\widehat{\lambda}_j : j\geq1\}$ are the population and empirical eigenvalues, $\{\phi_j(t):j\geq1 \}$ and $\{\widehat{\phi}_j(t):j\geq1 \}$ are the corresponding eigenfunctions, or functional principal components, each forming an orthonormal basis of $L_2([0,1])$. Then, $X_i$ and $\beta_k$
can be expanded as
\begin{equation}\label{}
  X_i(t)=\mu_X(t)+\sum_{j=1}^{\infty} \theta_{ij}\widehat{\phi}_j(t), \qquad \beta_k(t)=\sum_{j=1}^{\infty} b_{kj}\widehat{\phi}_j(t),
\end{equation}
where $\theta_{ij}=\langle X_i-\mu_X,\widehat{\phi}_j\rangle$ and $b_{kj}=\langle \beta_k,\widehat{\phi}_j\rangle$. We approximate $X_i$ and $\beta_k$ by
\begin{equation}\label{}
\begin{aligned}
  X_i^D(t)&=\mu_X(t)+\sum_{j=1}^{D} \theta_{ij}\widehat{\phi}_j(t)=\mu_X(t)+\bm{\theta}_i^T\bm{\Phi}_D(t), \\ \beta_k^D(t)&=\sum_{j=1}^{D} b_{kj}\widehat{\phi}_j(t)=\bm{b}_k^T\bm{\Phi}_D(t),
  \end{aligned}
\end{equation}
for some suitably large $D$, where $\bm{\theta}_i=(\theta_{i1},\ldots, \theta_{iD})^T$, $\bm{b}_k=(b_{k1},\ldots,b_{kD})^T$, $\bm{\Phi}_D(t)=(\widehat{\phi}_1(t),\ldots,\widehat{\phi}_D(t))^T$.

According to \cite{sheng2016sufficient},
the sample version of $\mathcal{V}^2(\langle \bm{\beta},X\rangle,Y)$ denoted by $\mathcal{V}_n^2(\langle \bm{\beta},\bm{X}\rangle,\bm{Y})$ has the following form:
\begin{equation}\label{}
  \mathcal{V}_n^2(\langle \bm{\beta},\bm{X}\rangle,\bm{Y})=\dfrac{1}{n^2}\sum_{i,j=1}^{n}A_{ij}(\bm{\beta})B_{ij},
\end{equation}
where
\begin{equation}\begin{array}{l}\label{sampleV2}
A_{ij}(\bm{\beta})=a_{ij}(\bm{\beta})-\bar{a}_{i .}(\bm{\beta})-\bar{a}_{. j}(\bm{\beta})+\bar{a}_{. .}(\bm{\beta}), \\\\
a_{ij}(\bm{\beta})=\left|\langle\bm{\beta}, X_{i}\rangle-\langle\bm{\beta}, X_{j}\rangle\right|_K,\quad \bar{a}_{i}(\bm{\beta})=\frac{1}{n} \sum_{j=1}^{n} a_{ij}(\bm{\beta}), \\\\
\bar{a}_{. j}(\bm{\beta})=\frac{1}{n} \sum_{i=1}^{n} a_{ij}(\bm{\beta}),\quad \bar{a}_{. .}(\bm{\beta})=\frac{1}{n^{2}} \sum_{i, j=1}^{n} a_{ij}(\bm{\beta}),
\end{array}\end{equation}
and $B_{ij}=b_{ij}-\bar{b}_{i.}-\bar{b}_{.j}+\bar{b}_{..}$, $b_{ij}=|Y_i-Y_j|$, and the definition of $\bar{b}_{i.}$, $\bar{b}_{.j}$, $\bar{b}_{..}$ is similar to those of $\bar{a}_{i.}(\bm{\beta})$, $\bar{a}_{.j}(\bm{\beta})$, $\bar{a}_{..}(\bm{\beta})$.
Then, an estimate of  a basis of $S_{y|x}$, say $\widehat{\bm{\eta}}_{n,D}$ is obtained by solving the following optimization problem:
\begin{equation}\label{opt}
  \widehat{\bm{\eta}}_{n,D}=\mathop{argmax}_{\langle\bm{\beta}^D ,\bm{\beta}^D\rangle_{\widehat{\Sigma}_X}=I_{K}} \mathcal{V}_n^2(\langle \bm{\beta}^D,\bm{X}^D\rangle,\bm{Y}),
\end{equation}
where $\bm{\beta}^D=\big(\beta_1^D(t),\ldots,\beta_K^D(t)\big)$ and $\bm{X}^D=\big(X_1^D,\ldots,X_n^D\big)^T$. We just need to replace $\langle\bm{\beta}, X_{i}\rangle$ and $\langle\bm{\beta}, X_{j}\rangle$ in \eqref{sampleV2} with $\langle\bm{\beta}^D, X_{i}^D\rangle$ and $\langle\bm{\beta}^D, X_{j}^D\rangle$ to calculate $\mathcal{V}_n^2(\langle \bm{\beta}^D,\bm{X}^D\rangle,\bm{Y})$. Particularly, let $\bm{B}=(\bm{b}_1,\ldots,\bm{b}_K)$, then
$\langle\bm{\beta}^D, X_{i}^D\rangle=\bm{B}^T\bm{\theta}_i$ and
$\langle\bm{\beta}^D , \bm{\beta}^D\rangle_{\widehat{\Sigma}_X}=\bm{B}^T\widehat{\Sigma}_{\bm{\theta}}\bm{B}$
where $\widehat{\Sigma}_{\bm{\theta}}=1/n\sum_{i=1}^{n}(\bm{\theta}_i-\bar{\bm{\theta}})(\bm{\theta}_i-\bar{\bm{\theta}})^T$,
$\bar{\bm{\theta}}=1/n\sum_{i=1}^{n}\bm{\theta}_i$. Thus, the optimization problem \eqref{opt} is equivalent to
\begin{equation}\label{maxB}
  \widehat{\bm{B}}_{n,D}=\mathop{argmax}_{\bm{B}^T\widehat{\Sigma}_{\bm{\theta}}\bm{B}=I_{K}} \mathcal{V}_n^2(\bm{B}^T\bm{\theta},\bm{Y}),
\end{equation}
where $\bm{\theta}=(\bm{\theta}_1,\ldots, \bm{\theta}_K)$ and $\widehat{\bm{B}}_{n,D}$ is a $D\times K$ matrix.
Similarly, we can replace $\langle\bm{\beta}, X_{i}\rangle$ and $\langle\bm{\beta}, X_{j}\rangle$ in \eqref{sampleV2} with
$\bm{B}^T\bm{\theta}_i$ and $\bm{B}^T\bm{\theta}_j$ to calculate  $\mathcal{V}_n^2(\bm{B}^T\bm{\theta},\bm{Y})$.
Although the optimization problem \eqref{maxi} with respect to $\beta_k$ is taken over an infinite dimensional space, the solution can actually be found in a finite dimensional subspace by regularization. It suffices to estimate the coefficients matrix $\bm{B}$ in \eqref{maxB}.

Note that it is complicated to solve the optimization problem \eqref{maxB} over a $D\times K$ matrix.
A projection pursuit type
of sufficient searching algorithm \cite{xue2017sufficient} is adopted to break down the problem \eqref{maxB} into
successive single-index searching.
The algorithm can be described as follows:
\begin{enumerate}
  \item[1.] Solve the single-index searching problem $\widehat{\gamma}_{1}=\arg\max_{\bm{B}^T\widehat{\Sigma}_{\bm{\theta}}\bm{B}=1} \mathcal{V}_n^2(\bm{B}^T\bm{\theta},\bm{Y})$, where $\bm{B}$ is a $D\times1$ vector. $\widehat{\gamma}_1$ is the first intermediate direction.
  \item[2.] Construct $D\times (D-1)$ matrix $\Gamma_1$ such that $\widehat{\Sigma}_{\bm{\theta}}^{1/2}(\widehat{\gamma}_1,\Gamma_1)$ is an orthogonal matrix.
  \item[3.] Let $a\in \mathbb{R}^{D-1}$ and consider the predictor matrix $(\widehat{\gamma}_{1},\Gamma_1a)$, where
  $\widehat{\gamma}_1$ is fixed. Solve the problem $$a_{1}=\underset{(\widehat{\gamma}_{1},\Gamma_1a)^T\widehat{\Sigma}_{\bm{\theta}}(\widehat{\gamma}_{1},\Gamma_1a)=I_2}{\arg\max}
  \left\{\mathcal{V}_n^2((\widehat{\gamma}_{1},\Gamma_1a)^T\bm{\theta},\bm{Y}):a\in \mathbb{R}^{D-1}\right\},$$
  then the second intermediate direction is $\widehat{\gamma}_2=\Gamma_1a_1$.
  \item[4.] Let the $D\times1$ vectors $\widehat{\gamma}_1,\widehat{\gamma}_2,\ldots,\widehat{\gamma}_k$ be the first $k$ intermediate directions, and let $\widehat{\Sigma}_{\bm{\theta}}^{1/2}(\widehat{\gamma}_1,\widehat{\gamma}_2,\ldots,\widehat{\gamma}_k,\Gamma_k)$ form an orthogonal matrix. Then we search for a $(D-k)\times1 $ vector $a_k$ based on the predictor matrix
      $(\widehat{\gamma}_1,\widehat{\gamma}_2,\ldots,\widehat{\gamma}_k,\Gamma_ka_k)$. Then the $(k+1)$-th intermediate direction is $\widehat{\gamma}_{k+1}=\Gamma_ka_k$.
  \item[5.] The estimate for the coefficients matrix $\bm{B}$ in \eqref{maxB} is
  $\widehat{\bm{B}}_{n,D}=(\widehat{\gamma}_1,\widehat{\gamma}_2,\ldots,\widehat{\gamma}_K)$.
\end{enumerate}

 Finally, the estimate of $S_{y|x}$ is
\begin{equation}\label{etand}
  \widehat{\bm{\eta}}_{n,D}=\bm{\Phi}_D^T(t)\widehat{\bm{B}}_{n,D}.
\end{equation}

For the asymptotic analysis,
it is typical to allow $D$ to diverge with the increase of $n$. We add a subscript $n$ to $D$ to
emphasize this relationship. The following theorem states the consistency of the estimator $\widehat{\bm{\eta}}_{n,D_n}$.

\begin{theorem}\label{thm1}
   Assume $\bm{\eta}=(\eta_1,\ldots,\eta_K)$ is a basis of $S_{y|x}$ with $\langle\bm{\eta} , \bm{\eta}\rangle_{\Sigma_X}=I_K$.
  For $\widehat{\bm{\eta}}_{n,D_n}$ defined in \eqref{opt}, as $n\rightarrow\infty$, $D_n\rightarrow\infty$ we have $\widehat{\bm{\eta}}_{n,D_n} \stackrel{P}{\longrightarrow} \bm{\eta}$, provided that Condition \ref{assump1} holds.
\end{theorem}

Theorem \ref{thm1} establishes the consistency of the estimator $\widehat{\bm{\eta}}_{n,D_n}$ of $\bm{\eta}$.
Here, we denote $\widehat{\bm{\eta}}_{n,D_n}=(\widehat{\eta}_{n,D_n,1},\ldots,\widehat{\eta}_{n,D_n,K})$, and the expression
$\widehat{\bm{\eta}}_{n,D_n} \stackrel{P}{\longrightarrow} \bm{\eta}$ means $\|\widehat{\bm{\eta}}_{n,D_n}-\bm{\eta}\|=\big(\sum_{k=1}^{K}\|\widehat{\eta}_{n,D_n,k}-\eta_k\|^2\big)^{1/2}\stackrel{P}{\longrightarrow}0$.

\subsection{Estimation for spase longitudinal data}
The focus of this subsection is to estimate $S_{y|x}$ for intermittently and sparsely measured longitudinal covariates.
When only
a few observations are available for some or even all subjects, individual smoothing to recover
$X_i$ is infeasible and one must pool data across subjects for consistent estimation. For the i.i.d sample $(\bm{X},\bm{Y})=\{(X_i(t),Y_i),\ i=1,2,\ldots,n\}$, the predictors $X_i$ are observed intermittently, contaminated with noise, and observed in the form of $\big\{(T_{i j}, U_{i j}): i=1, \ldots, n ; j=1, \ldots, N_{i}\big\}$ where
\begin{equation*}\label{}
  U_{ij}=X_i(T_{ij})+\varepsilon_{ij}.
\end{equation*}
The i.i.d. measurement error $\varepsilon_{ij}$ satisfies $E(\varepsilon_{ij})=0$ and $var(\varepsilon_{ij})=\sigma^2_{\varepsilon}$. The numbers of observations $\{N_i\}_{i=1}^n$ are assumed to be random,
reflecting sparse and irregular designs. The observation time points $\{T_{ij}\}$ are assumed to be i.i.d. realizations of
a random variable and independent of all other random variables. Another assumption is that the pooled time points
$\{T_{ij }\}$ are sufficiently dense in the domain of $X(t)$.

Firstly, we estimate the mean function $\mu_X(t)$ based on the pooled data. Following \cite{yao2005functional}, local linear smoothing is conducted for estimating $\mu_X(t)$ by minimizing
\begin{equation}\label{}
  \sum_{i=1}^{n}\sum_{j=1}^{N_i} K_1\left(\frac{T_{ij}-t}{h_{\mu}}\right)\left\{U_{ij}-a_0-a_1(t-T_{ij})\right\}^2,
\end{equation}
with respect to $a_0$ and $a_1$, where $K_1$ is a univariate kernel function and $h_{\mu}$ is the bandwidth. Then the estimate of $\mu_X(t)$ is $\widehat{\mu}_X(t)=\widehat{a}_0$.
For the covariance function $\Sigma_X(s,t)$, \cite{yao2005functional} defined the observed raw covariance by
$G_i(T_{ij},T_{il})=(U_{ij}-\widehat{\mu}_X(T_{ij}))(U_{il}-\widehat{\mu}_X(T_{il}))$. Solving the local linear surface
smoothing problem
\begin{equation}
\min _{\left(b_{0}, b_{1}, b_{2}\right)} \sum_{i=1}^{n} \sum_{j \neq l}^{N_{i}}K_{2}\left(\frac{T_{i j}-s}{h_{\Sigma}}, \frac{T_{i l}-t}{h_{\Sigma}}\right)\left\{G_{i}\left(T_{i j}, T_{i l}\right)-b_{0}-b_{1}\left(T_{i j}-s\right)-b_{2}\left(T_{i l}-t\right)\right\}^{2},
\end{equation}
yields $\widehat{\Sigma}_X(s,t)=\widehat{b}_0$, where $K_2$ is a bivariate kernel function with bandwidth $h_{\Sigma}$.
Then, similar to the previous subsection, $\widehat{\lambda}_j$ and $\widehat{\phi}_j$ can be obtained from the eigen decomposition of $\widehat{\Sigma}_X(s,t)$.

From the optimization problem \eqref{maxB}, we known that the only quantity we need to acquire is the FPC scores
$\theta_{ij}=\langle X_i-\mu_X,{\phi}_j\rangle$. For sparse longitudinal data, individual smoothing to recover
$X_i$ is infeasible, thus numerical integration for calculating $\theta_{ij}$ will not provide reasonable approximations to
the real FPC scores. We adopt the efficient \emph{Principal Analysis by Conditional Expectation (PACE)} \citep{yao2005functional} method specifically
designed for sparse longitudinal data to estimate FPC scores. Denote $\widetilde{\bm{X}}_i=(X_i(T_{i1}),\ldots,X_i(T_{iN_i}))^T$, $\widetilde{\bm{U}}_i=(U_{i1},\cdots,U_{iN_i})^T$,
$\bm{\mu}_i=(\mu_X(T_{i1}),\ldots,\mu_X(T_{iN_i}))^T$, and $\bm{\phi}_{ij}=(\phi_j(T_{i1}),\ldots,\phi_j(T_{iN_i}))^T$.
When $\theta_{ij}$ and $\varepsilon_{ij}$ are jointly Gaussian, the best prediction of the FPC score $\theta_{ij}$
given $i$-th subject is the conditional expectation
\begin{equation}\label{}
  \widetilde{\theta}_{ij}=E(\theta_{ij}|\widetilde{\bm{U}}_i)=
  \lambda_j\bm{\phi}_{ij}^T\bm{\Sigma_{U_i}}^{-1}(\widetilde{\bm{U}}_i-\bm{\mu}_i),
\end{equation}
where $\bm{\Sigma_{U_i}}=cov(\widetilde{\bm{U}}_i,\widetilde{\bm{U}}_i)
=cov(\widetilde{\bm{X}}_i,\widetilde{\bm{X}}_i)+\sigma_{\varepsilon}^2\bm{I}_{N_i}$ and the $N_i\times N_i$ matrix $cov(\widetilde{\bm{X}}_i,\widetilde{\bm{X}}_i)$ =
$\Big(\Sigma_X(T_{ij},T_{il})\Big)_{1\leq j,l \leq N_i}$. By substituting estimates of $\bm{\mu}_i$, $\lambda_j$, $\bm{\phi}_{ij}$ and $\bm{\Sigma_{U_i}}$ obtained from the pooled data, we have an estimate of $\theta_{ij}$,
\begin{equation}\label{}
  \widehat{\theta}_{ij}=\widehat{E}(\theta_{ij}|\widetilde{\bm{U}}_i)=
  \widehat{\lambda}_j\widehat{\bm{\phi}}_{ij}^T\widehat{\bm{\Sigma}}_{\bm{U}_i}^{-1}(\widetilde{\bm{U}}_i-\widehat{\bm{\mu}}_i),
\end{equation}
where $\widehat{\bm{\Sigma}}_{U_i}=\widehat{cov}(\widetilde{\bm{X}}_i,\widetilde{\bm{X}}_i)+\widehat{\sigma}_{\varepsilon}^2\bm{I}_{N_i}$,
 $\widehat{cov}(\widetilde{\bm{X}}_i,\widetilde{\bm{X}}_i)=\Big(\widehat{\Sigma}_X(T_{ij},T_{il})\Big)_{1\leq j,l \leq N_i}$, and $\widehat{\sigma}_{\varepsilon}^2$ is an estimate of ${\sigma}_{\varepsilon}^2$ which can be found in \cite{yao2005functional}. In this paper, we do not introduce the specific estimate of ${\sigma}_{\varepsilon}^2$ for the sake of brevity. \cite{yao2005functional} has shown that under some regularization conditions,
$\widehat{\theta}_{ij}\stackrel{P}{\longrightarrow}\widetilde{\theta}_{ij}$, which means that $\widehat{\theta}_{ij}$ is a good estimate of $\theta_{ij}$.
We substitute $\{\widehat{\theta}_{ij}\}$ into the optimization problem \eqref{maxB} to get the coefficients matrix $\widehat{\bm{B}}_{n,D}$. The optimization steps are exactly the same as those of completely observed  case.

\subsection{Selection of tuning parameters}\label{subsecK}

Computation of  the proposed method relies on the choice of two parameters: the truncation number $D$ and the
structural dimension $K$.
Determining $D$ and $K$ can be tackled in different ways and it depends on the goal
of the analysis. If the purpose is prediction, $D$ and $K$ can be treated
as parameters of the whole model and adjusted according to the performance
of the prediction, such as cross-validation.
This has been successfully experimented with in applications \citep{ferre2005smoothed,hsing2009rkhs,yao2015effective}.
Consider  model \eqref{fSDR}:
\begin{equation*}
  Y=g\left(\left\langle\beta_{1}, X\right\rangle, \ldots,\left\langle\beta_{K}, X\right\rangle, \epsilon\right).
\end{equation*}
The data $(\bm{X},\bm{Y})$ were randomly divided into $k$ equal portions
$\{(\bm{X}^{(1)},\bm{Y}^{(1)}),\ldots,(\bm{X}^{(k)},\bm{Y}^{(k)})\}$.
For each feasible $D$, $K$ and $i=1,\ldots,k$, we leave out $(\bm{X}^{(i)},\bm{Y}^{(i)})$
and use the rest of the data $(\bm{X}^{(-i)},\bm{Y}^{(-i)})$ to compute the
$\widehat{\bm{\eta}}_{n,D}$ in \eqref{etand} and nonparametrically estimate $g$.
Then, use the $\widehat{\bm{\eta}}_{n,D}$, the estimate $g$, and $\bm{X}^{(i)}$ to compute
predicted values $\widehat{\bm{Y}}_{D, K}^{(i)}$.
Let $\mathrm{CV}(D,K)=1/k\cdot\sum_{i=1}^{k}\left\|\bm{Y}^{(i)}-\widehat{\bm{Y}}_{D, K}^{(i)}\right\|_2^2$,
where $\|\cdot\|_2$ is the Euclidean norm, and pick $D$ and $K$ to minimize $\mathrm{CV}(D,K)$.
However, the cross-validation procedures are not ideal since the nonparametric fitting adds an extra layer of complication.
A more satisfactory of the selection of $D$ and $K$ is currently not available.
We will apply this  cross-validation method  in the real data analysis.

When  FSDR is used in  a  descriptive way,  we are mainly interested in recovering the directions per se.
Much of the existing literature, such as \cite{hsing2009rkhs, wang2014functional, lian2015functional},
 suggested  that $D$ can be chosen subjectively. Specifically,
the number $D$ of included eigenfunctions is chosen by fraction of variance explained criterion in practice,
$$D=\min \left\{k: \sum_{l=1}^{k} \hat{\lambda}_{l} \Big/ \sum_{l=1}^{n} \hat{\lambda}_{l} \geq R\right\},$$
with a given threshold $R$ close to $1$ and the eigenvalues $\lambda_l, \ 1\leq l\leq k$ are not ``too small".
We can also draw a scree plot and choose the ``elbow" point as the truncation number.
In the simulations in this paper,
we recommend $R=95\%$, which includes $5$ eigenfunctions, and this choice of $D$ yields satisfactory results.

We then consider the selection of $K$ when $D$ is known.
Similar as in the multivariate case, the selection of $K$ relies on a criterion measuring the quality
of the estimation of $S_{y|x}$.
A bootstrap method suggested by \cite{sheng2016sufficient} and \cite{zhang2019robust} can be readily extended to our method and we apply it in the real data analysis.
Specifically, the bootstrap method is based on the measure of distance of two functional spaces \citep{yao2015effective}:
 \begin{equation}\label{}
   \Delta_m(S_1,S_2)=\|P_{S_1}-P_{S_2}\|_H,
 \end{equation}
where $S_i$ is spanned by $\{\beta_1,\ldots,\beta_{d_i}\}$, $P_{S_i}=\sum_{j=1}^{d_i}\beta_j\otimes\beta_j$ is the projection operator onto $S_i$ for
$i=1,2$, and $\|A\|^2_H=\int\int A^2(s,t)ds\ dt$ for a linear operator $A\in L_2([0,1]\times[0,1])$. Obviously, the smaller the distance is, the closer the two spaces are.

In order to use this measure, we will treat $(\bm{\theta},\bm{Y})=\{(\bm{\theta}_i,Y_i),i=1,\ldots,n\}$ as the sample. For each possible working dimensions $1\leq k\leq D-1$, we solve the problem \eqref{maxB} to obtain an estimated coefficient matrix $\widehat{\bm{B}}_k$ whose columns spans a subspace in $\mathbb{R}^D$ and then we get the estimator of $S_{y|x}$, $\widehat{\bm{\eta}}_{k}=\bm{\Phi}_D^T(t)\widehat{\bm{B}}_{k}$. Here we omit the subscript $n$ and $D$ to emphasize the status of $k$.
Then we randomly sample the data $(\bm{\theta},\bm{Y})$ with replacement $B$ times, and obtain the estimated subspace based on the bootstrap samples, and denote them by $\widehat{\bm{\eta}}^b_k,\ b=1,\ldots,B$. We calculate
$\Delta_m(\widehat{\bm{\eta}}_k,\widehat{\bm{\eta}}_k^b)$ for $b=1,\ldots,B$ and use the mean $1/B\cdot\sum_{b=1}^{B}\Delta_m(\widehat{\bm{\eta}}_k,\widehat{\bm{\eta}}_k^b)$ as the measure of variability for each $k$.
We choose the $k$ corresponding to the smallest variability as our estimated $K$. The reason why this bootstrap method works well is mentioned in \cite{sheng2016sufficient}.

It is worth noting that directly using the bootstrap method  is  time consuming.
However, the bootstrap method can easily be modified to a parallel version to significantly reduce the computation
time. We use the \texttt{R} package \texttt{parallel} for parallel computing.

\begin{figure}[b]
\centering
\subfigure{
\includegraphics[width=3.5cm]{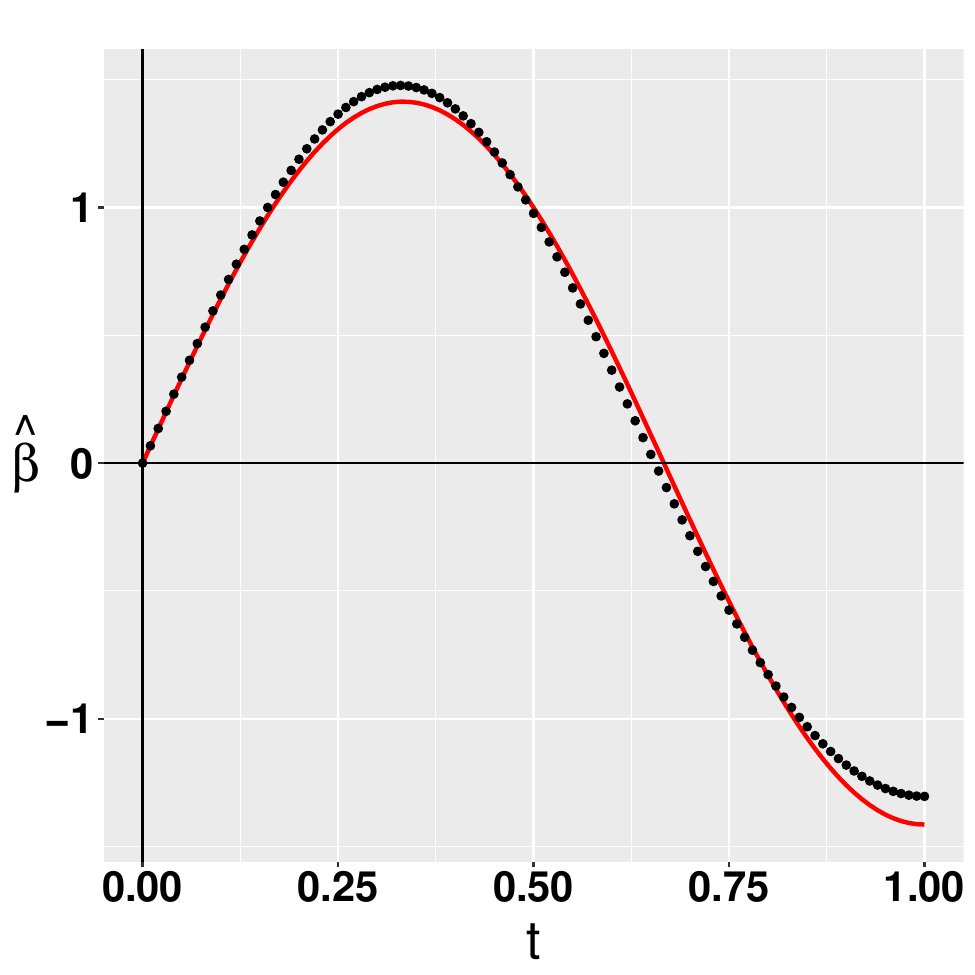}
}
\quad
\subfigure{
\includegraphics[width=3.5cm]{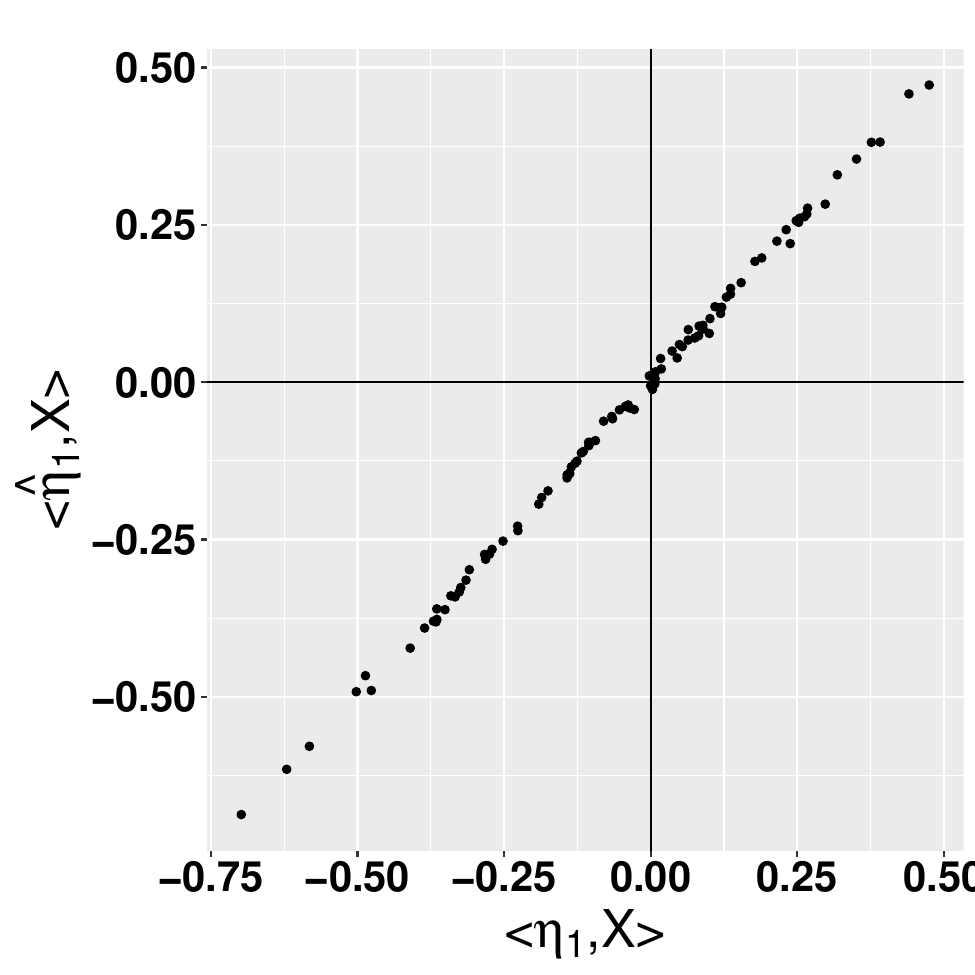}
}
\quad
\subfigure{
\includegraphics[width=3.5cm]{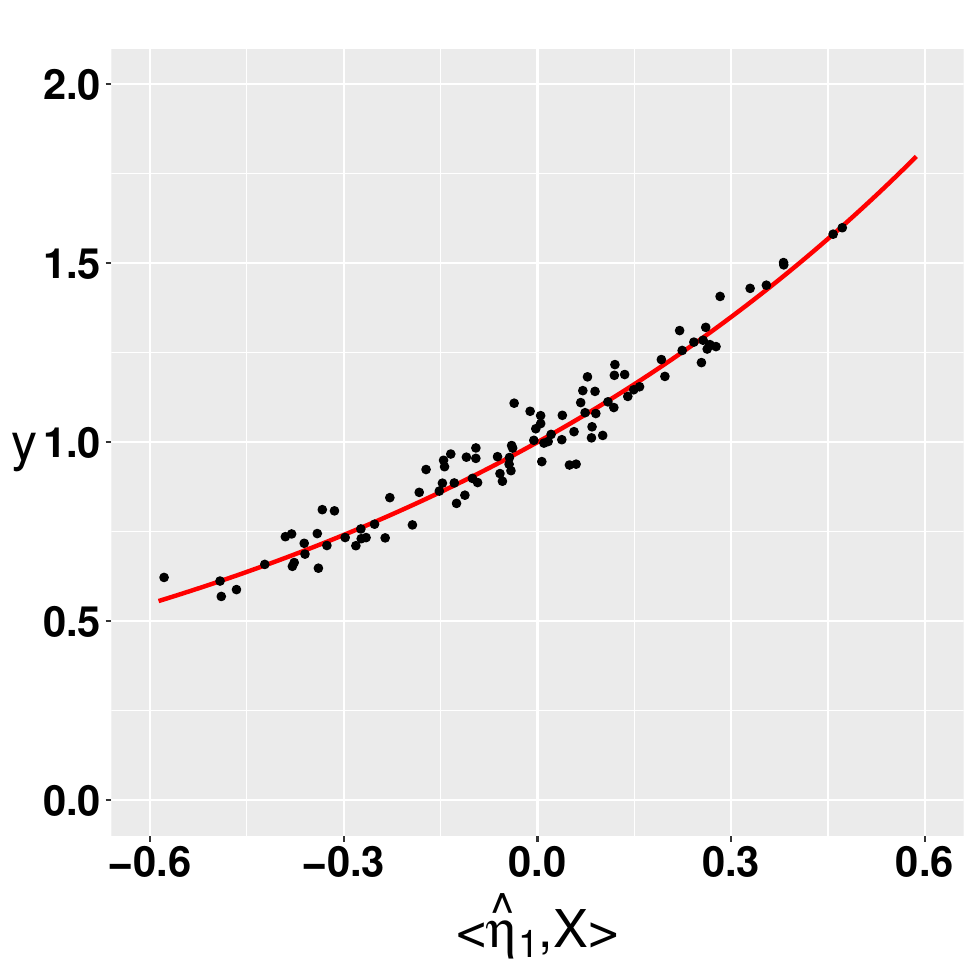}
}
\caption{Left panel: $\eta_1$ (red smooth curve) and $\widehat{\eta}_1$ (black dots) versus $t$.
Middle panel: $\langle \widehat{\eta}_1, X\rangle$ versus $\langle \eta_1, X\rangle$.
Right panel: $Y$ versus $\langle \widehat{\eta}_1, X\rangle$ and the true link function $y=e^x$
(red smooth curve) .}\label{model1fig}
\end{figure}

\section{Numerical studies}\label{numerical}
\subsection{Simulations}
In this subsection, we conduct simulation studies to provide an insight in the empirical performance of the proposed  method (FDCOV) and compare it with some existing methods. For the case of  completely observed functional data, the competitors
includes functional sliced inverse regression
(FSIR, \citealp{ferre2003functional}), functional sliced average variance estimation (FSAVE, \citealp{lian2014series}),
functional contour regression (FCR, \citealp{wang2013functional}) and
functional directional regression (FDR, \citealp{wang2020directional}).
For the case of  sparse longitudinal data, we compare our proposed  method with functional inverse regression (FIR, \citealp{jiang2014inverse}) and
functional cumulative slicing (FCS, \citealp{yao2015effective}).

In the simulation studies, the following five models are considered:
\begin{equation*}\begin{aligned}
&\text {(1) } Y=\exp(\left\langle\eta_{1}, X\right\rangle)+\epsilon,\quad \eta_{1}(t)=\sin (3 \pi t / 2), \\
&\text {(2) } Y=\exp(\left\langle\eta_{1}, X\right\rangle)+\exp(|\left\langle\eta_{2}, X\right\rangle|)+\epsilon,\ \eta_{1}(t)=\sin (3 \pi t / 2),\ \eta_{2}(t)=\sin (5 \pi t / 2) , \\
&\text {(3) } Y=\exp(\left\langle\eta_{1}, X\right\rangle)+\exp(|\left\langle\eta_{2}, X\right\rangle|)+\epsilon,\ \eta_{1}(t)=(2 t-1)^{3}+1,\ \eta_{2}(t)=\cos ((2t-1) \pi )+1 ,\\
&\text {(4) } Y=5 \left\langle \eta_{1},X\right\rangle +15\left\langle \eta_{2},X\right\rangle^{2}\epsilon,\ \eta_{1}(t)=\sin (3 \pi t/2),\ \eta_{2}(t)=\sin (5 \pi t/2) ,\\
&\text {(5) } Y=50\left\langle\eta_{1}, X\right\rangle\left\langle\eta_{2}, X\right\rangle^{2}+\epsilon,\ \eta_{1}(t)=(2 t-1)^{2}-1,\ \eta_{2}(t)=\sin (5 \pi t / 2),
\end{aligned}\end{equation*}
where $\epsilon \sim N(0,0.1^2)$ and $X$ is the standard Brownian motion on $[0,1]$, independent of $\epsilon$. These examples cover various situations. Model 1 is a single index model from \cite{hsing2009rkhs}. Model 2 and Model 3 are taken from \cite{ferre2003functional}. In Model 2, both $\eta_1$ and $\eta_2$ are eigenvectors of the Brownian motion, while both $\eta_1$ and $\eta_2$ in Model 3 are not eigenvectors of the Brownian motion. Model 4 was used in \cite{lian2015functional}, considering  heterogeneous errors. Model 5
was previously considered in \cite{lian2014series} and \cite{lian2015functional}, where $\eta_1$ is not a eigenvector  of the Brownian motion and the link function is not additive.

\begin{table}[t]
\centering
\begin{threeparttable}[h]
\caption{Estimation error $\|\widehat{P}-P\|_H$ for different estimators for Model 1. Numbers in parentheses are standard errors calculated from
100 generated data sets.}\label{m1}
\begin{tabular}{cccccccc}
\hline
$n$ & $L$& $r$  & FDCOV        & FSIR         & FSAVE       & FDR         & FCR      \\ \hline
100 & 5  & 0.05 & 0.121(0.047) & 0.255(0.117) & 0.287(0.128)& 0.253(0.129)& 0.174(0.076) \\
    & 10 & 0.10 & 0.121(0.047) & 0.190(0.077) & 0.379(0.243)& 0.296(0.127)& 0.168(0.073)    \\
    & 15 & 0.15 & 0.121(0.047) & 0.169(0.072) & 1.127(0.400)& 0.359(0.162)& 0.160(0.069)    \\
    & 20 & 0.20 & 0.121(0.047) & 0.177(0.077) & 1.299(0.248)& 0.415(0.162)& 0.163(0.073)    \\
200 & 5  & 0.05 & 0.085(0.036) & 0.165(0.060) & 0.166(0.066)& 0.164(0.071)& 0.120(0.050)    \\
    & 10 & 0.10 & 0.085(0.036) & 0.130(0.062) & 0.154(0.064)& 0.194(0.094)& 0.115(0.053)    \\
    & 15 & 0.15 & 0.085(0.036) & 0.114(0.054) & 0.195(0.149)& 0.201(0.084)& 0.113(0.050)    \\
    & 20 & 0.20 & 0.085(0.036) & 0.123(0.050) & 0.239(0.190)& 0.232(0.114)& 0.116(0.048)    \\ \hline
\end{tabular}
\end{threeparttable}
\end{table}

\begin{table}[t]
\centering
\begin{threeparttable}[h]
\caption{Estimation error $\|\widehat{P}-P\|_H$ for different estimators for Model 2. Numbers in parentheses are standard errors calculated from
100 generated data sets.}\label{m2}
\begin{tabular}{cccccccc}
\hline
$n$ & $L$& $r$  & FDCOV        & FSIR         & FSAVE       & FDR         & FCR      \\ \hline
100 & 5  & 0.05 & 0.458(0.112) & 1.414(0.273) & 1.208(0.297)& 0.970(0.255)& 1.069(0.312) \\
    & 10 & 0.10 & 0.458(0.112) & 1.444(0.277) & 1.359(0.216)& 1.044(0.263)& 0.992(0.275)    \\
    & 15 & 0.15 & 0.458(0.112) & 1.482(0.278) & 1.453(0.180)& 1.334(0.276)& 0.921(0.246)    \\
    & 20 & 0.20 & 0.458(0.112) & 1.512(0.297) & 1.497(0.209)& 1.288(0.240)& 0.916(0.264)    \\
200 & 5  & 0.05 & 0.187(0.075) & 1.369(0.255) & 0.874(0.314)& 0.711(0.270)& 0.756(0.181)    \\
    & 10 & 0.10 & 0.187(0.075) & 1.294(0.298) & 1.059(0.304)& 0.693(0.231)& 0.724(0.182)    \\
    & 15 & 0.15 & 0.187(0.075) & 1.316(0.305) & 1.184(0.305)& 0.793(0.242)& 0.718(0.185)    \\
    & 20 & 0.20 & 0.187(0.075) & 1.392(0.270) & 1.341(0.213)& 0.870(0.267)& 0.687(0.194)    \\ \hline
\end{tabular}
\end{threeparttable}
\end{table}

\begin{table}[t]
\centering
\begin{threeparttable}[h]
\caption{Estimation error $\|\widehat{P}-P\|_H$ for different estimators for Model 3. Numbers in parentheses are standard errors calculated from
100 generated data sets.}\label{m3}
\begin{tabular}{cccccccc}
\hline
$n$ & $L$& $r$  & FDCOV        & FSIR         & FSAVE       & FDR         & FCR      \\ \hline
100 & 5  & 0.05 & 1.733(0.028) & 2.712(0.102) & 2.288(0.251)& 2.274(0.265)& 2.101(0.228) \\
    & 10 & 0.10 & 1.733(0.028) & 2.698(0.128) & 2.444(0.255)& 2.396(0.269)& 2.122(0.229)    \\
    & 15 & 0.15 & 1.733(0.028) & 2.729(0.120) & 2.499(0.236)& 2.481(0.275)& 2.121(0.206)    \\
    & 20 & 0.20 & 1.733(0.028) & 2.713(0.123) & 2.569(0.272)& 2.497(0.228)& 2.180(0.225)    \\
200 & 5  & 0.05 & 1.724(0.017) & 2.710(0.094) & 2.126(0.217)& 2.107(0.205)& 1.939(0.147)    \\
    & 10 & 0.10 & 1.724(0.017) & 2.678(0.131) & 2.193(0.235)& 2.166(0.227)& 1.974(0.156)    \\
    & 15 & 0.15 & 1.724(0.017) & 2.683(0.125) & 2.317(0.236)& 2.162(0.206)& 2.004(0.190)    \\
    & 20 & 0.20 & 1.724(0.017) & 2.692(0.123) & 2.268(0.255)& 2.219(0.231)& 2.020(0.164)    \\ \hline
\end{tabular}
\end{threeparttable}
\end{table}
For the considered methods for completely observed data, we need to decide the following tuning parameters.
For FSIR, FSAVE and FDR: the number of slices $L$; for FCR: the proportion $r$ of empirical directions. We consider $L=5,\ 10,\ 15,\ 20$, $r=0.05,\ 0.10,\ 0.15,\ 0.20$ and $n=100,\ 200$. In each setting, we simulate 100 data sets and each
random curve is sampled at $p = 100$ equally spaced points in $[0, 1]$, $\{t_1,\ldots,t_{100}\}$ with $t_1=0$ and $t_{100}=1$.
Similar to the setting in \cite{wang2020directional},
all these methods are implemented using functional principal component analysis with truncation $D$ chosen such that
$95\%$ of variability in the predictor are retained. That is
\begin{equation}D=\min\left\{k:\left(\sum_{i=1}^{k} \widehat{\lambda}_{i}\right) \Bigg/\left(\sum_{i=1}^{n} \widehat{\lambda}_{i}\right) \geq 95 \%\right\}.\end{equation}

We assume that the structural dimension $K$ is known. Let $P=\sum_{k=1}^{K}\eta_k\otimes\eta_k$ and
$\widehat{P}=\sum_{k=1}^{K}\eta_{n,D,k}\otimes\eta_{n,D,k}$ be the projection operators onto the true $S_{y|x}$ and estimated $S_{y|x}$ respectively. We calculate $\|\widehat{P}-P\|_H$ as the estimation error with smaller values indicating better estimation performance. All the simulation results for completely observed data are reported in Figure \ref{model1fig} and Tables \ref{m1}-\ref{m5}.

For Model 1, some results of our proposed method are displayed in Figure \ref{model1fig}.
The left panel of Figure \ref{model1fig} gives the plots of $\eta_1$ and $\widehat{\eta}_1$. The red smooth curve is the true direction $\eta_1=\sin(3\pi t/2)$ and the black dots are corresponding estimator $\widehat{\eta}_1$ at $100$ equally spaced time points in $[0,1]$. We see that $\widehat{\eta}_1$ coincides almost perfectly with $\eta_1$. In prediction stage, what we care about is that the estimated projection $\langle \widehat{\eta}_1, X\rangle$ is as close as possible to the true projection
$\langle \eta_1, X\rangle$. We plot in the middle panel of Figure \ref{model1fig} the indexes $\langle \widehat{\eta}_1, X\rangle$ versus $\langle \eta_1, X\rangle$. We find that these scatter plots reveal a strong correlation between
 $\langle \widehat{\eta}_1, X\rangle$ and $\langle \eta_1, X\rangle$. We also present the plot for $Y$ versus  $\langle \widehat{\eta}_1, X\rangle$ along with the true link function $y=e^x$ in the right panel of Figure \ref{model1fig}.

\begin{table}[t]
\centering
\begin{threeparttable}[]
\caption{Estimation error $\|\widehat{P}-P\|_H$ for different estimators for Model 4. Numbers in parentheses are standard errors calculated from
100 generated data sets.}\label{m4}
\begin{tabular}{cccccccc}
\hline
$n$ & $L$& $r$  & FDCOV        & FSIR         & FSAVE       & FDR         & FCR      \\ \hline
100 & 5  & 0.05 & 0.606(0.124) & 1.464(0.294) & 1.214(0.282)& 1.049(0.268)& 0.884(0.239) \\
    & 10 & 0.10 & 0.606(0.124) & 1.505(0.258) & 1.302(0.305)& 1.081(0.317)& 0.806(0.236)    \\
    & 15 & 0.15 & 0.606(0.124) & 1.493(0.316) & 1.454(0.247)& 1.460(0.320)& 0.830(0.219)    \\
    & 20 & 0.20 & 0.606(0.124) & 1.493(0.329) & 1.509(0.220)& 1.314(0.282)& 0.826(0.205)    \\
200 & 5  & 0.05 & 0.356(0.108) & 1.415(0.210) & 0.865(0.279)& 0.769(0.220)& 0.695(0.215)    \\
    & 10 & 0.10 & 0.356(0.108) & 1.383(0.234) & 0.965(0.286)& 0.805(0.226)& 0.674(0.156)    \\
    & 15 & 0.15 & 0.356(0.108) & 1.371(0.245) & 1.090(0.280)& 0.803(0.200)& 0.638(0.157)    \\
    & 20 & 0.20 & 0.356(0.108) & 1.383(0.265) & 1.135(0.293)& 0.809(0.255)& 0.659(0.181)    \\ \hline
\end{tabular}
\end{threeparttable}
\end{table}

\begin{table}[t]
\centering
\begin{threeparttable}[]
\caption{Estimation error $\|\widehat{P}-P\|_H$ for different estimators for Model 5. Numbers in parentheses are standard errors calculated from
100 generated data sets.}\label{m5}
\begin{tabular}{cccccccc}
\hline
$n$ & $L$& $r$  & FDCOV        & FSIR         & FSAVE       & FDR         & FCR      \\ \hline
100 & 5  & 0.05 & 0.632(0.052) & 1.147(0.063) & 1.041(0.202)& 1.022(0.133)& 1.037(0.227) \\
    & 10 & 0.10 & 0.632(0.052) & 1.140(0.051) & 1.314(0.362)& 1.044(0.131)& 1.006(0.208)    \\
    & 15 & 0.15 & 0.632(0.052) & 1.146(0.062) & 1.767(0.283)& 1.146(0.188)& 1.029(0.185)    \\
    & 20 & 0.20 & 0.632(0.052) & 1.146(0.065) & 1.807(0.250)& 1.246(0.311)& 1.039(0.207)    \\
200 & 5  & 0.05 & 0.605(0.035) & 1.121(0.036) & 0.813(0.128)& 0.930(0.113)& 0.841(0.163)    \\
    & 10 & 0.10 & 0.605(0.035) & 1.122(0.050) & 0.882(0.137)& 0.967(0.105)& 0.832(0.156)    \\
    & 15 & 0.15 & 0.605(0.035) & 1.122(0.059) & 1.018(0.236)& 0.973(0.111)& 0.887(0.195)    \\
    & 20 & 0.20 & 0.605(0.035) & 1.114(0.055) & 1.115(0.277)& 0.979(0.108)& 0.909(0.148)    \\ \hline
\end{tabular}
\end{threeparttable}
\end{table}

 Note that the results in Figure \ref{model1fig} are based on one single simulation run.  In order to get more representative results,
we compare our proposed method with FSIR, FSAVE, FDR and FCR based on 100 Monte Carlo repetitions. Tables \ref{m1}--\ref{m5}
report the mean and standard errors of $\|\widehat{P}-P\|_H$ for Models 1--5. From Tables \ref{m1}--\ref{m5}, FDCOV has the best performance in all five cases. For Model 1, FSAVE does not work well since it is known that SAVE is not efficient  in estimating monotone trends for small to moderate data sets. For Model 2 and Model 3, it is not surprising that the absolute value of $\|\widehat{P}-P\|_H$ of Model 2 is greater than that of Model 3 since Model 2 corresponding the ideal situation where
the true direction is included into the \emph{a priori} projection subspace. From Tables \ref{m2} and \ref{m4}, we see that
FDCOV can identify $S_{y|x}$ in heteroscedastic models, but it is not as efficient as in homoscedastic models.
For all cases, the results become better as
$n$ increases and the results are generally not very sensitive to the choice of $L$ and $r$. Note that FDCOV
has no parameters to tune and is not related to $L$ and $r$.

To generate the sparse longitudinal data, we randomly selected 10 to 20 observations from $\{t_1,t_2,...,t_{100}\}$
for each sample trajectory.  The measurement error $\varepsilon_{ij}$ is independent and identically distributed as
$N(0, 0.1^2)$.  The simulation consists of 100 runs and Table \ref{sparse} summarizes the results when $n$
is 100 and 200. For comparison, we also include the results of FIR \citep{jiang2014inverse} and FCS \citep{yao2015effective}.
The estimation of functional principal components for spare longitudinal data is implemented through \texttt{fdapace}
package in $\texttt{R}$ system.
The results suggest that our proposed method slightly
outperforms the other methods we considered.

\begin{table}[t]
\centering
\caption{Estimation error $\|\widehat{P}-P\|_H$ for different estimators. Numbers in parentheses are standard errors calculated from
100 generated data sets.}\label{sparse}
\begin{tabular}{ccccccc}
\hline
$n  $                   & methods& Model 1      & Model 2      &Model 3      & Model 4     & Model 5 \\ \hline
\multicolumn{1}{c}
{\multirow{4}{*}{$100$}}& FDCOV  &0.465(0.056)&0.752(0.060)&1.874(0.140)&0.986(0.211)&0.917(0.158)        \\
\multicolumn{1}{c}{}    & FIR    &1.007(0.060)&2.016(0.087)&2.532(0.697)&1.985(0.326)&1.282(0.265)        \\
\multicolumn{1}{c}{}    & FCS    &0.993(0.055)&2.004(0.062)&2.517(0.401)&1.998(0.326)&1.300(0.276)        \\ \hline
\multirow{4}{*}{$200$}  & FDCOV  &0.416(0.055)&0.718(0.054)&1.688(0.119)&0.851(0.177)&0.834(0.102)        \\
                        & FIR    &0.996(0.061)&1.991(0.070)&2.518(0.625)&1.887(0.324)&1.183(0.213)        \\
                        & FCS    &0.984(0.049)&1.996(0.059)&2.507(0.400)&1.897(0.245)&1.296(0.190)        \\\hline
\end{tabular}
\end{table}

\begin{table}[t]
\centering
\begin{threeparttable}[t]
\caption{Average distances using bootstrap samples for Models 1-5 }\label{KKK}
\begin{tabular}{cclllll}
\hline
Model              & $B$   & $K=1$  & $K=2$ & $K=3$     & $K=4$  & $K=5$                  \\
\hline
\multirow{2}{*}{$1$} & 100 & $0.083^*$ & $0.216$   & $0.404$ & $0.473$ & $0.517$ \\
                       & 200 & $0.076^*$ & $0.228$   & $0.398$ & $0.507$ & $0.545$ \\
\multirow{2}{*}{$2$} & 100 & $0.342$   & $0.106^*$ & $0.449$ & $0.498$ & $0.533$ \\
                       & 200 & $0.343$   & $0.097^*$ & $0.431$ & $0.472$ & $0.565$  \\
\multirow{2}{*}{$3$} & 100 & $1.838$   & $0.833^*$ & $1.609$ & $1.718$ & $1.390$   \\
                       & 200 & $1.690$   & $0.731^*$ & $1.439$ & $1.298$ & $1.092$ \\
\multirow{2}{*}{$4$} & 100 & $0.321$   & $0.117^*$ & $0.421$ & $0.550$ & $0.560$\\
                       & 200 & $0.389$   & $0.099^*$ & $0.481$ & $0.502$ & $0.528$ \\
\multirow{2}{*}{$5$} & 100 & $0.411$   & $0.198^*$ & $0.470$ & $0.598$ & $0.618$  \\
                       & 200 & $0.418$   & $0.201^*$ & $0.417$ & $0.605$ & $0.611$ \\
\hline
\end{tabular}
\begin{tablenotes}
     \item[] \footnotesize NOTE: A value with $*$ means it is the minimum average distance, which also corresponds to the selected
dimension.
   \end{tablenotes}
\end{threeparttable}
\end{table}

To examine the effectiveness of the bootstrap method for estimating $K$,  we still consider
the above five models. For each model, we consider $n=200$ and $B=100,200$.
As mentioned in Subsection \ref{subsecK}, we use the average distance $1/B\cdot\sum_{b=1}^{B}\Delta_m(\widehat{\bm{\eta}}_k,\widehat{\bm{\eta}}_k^b)$ as the measure of variability for each candidate $K$ and
the results of these average distances under different settings are summarized in Table \ref{KKK}.
The results show that the bootstrap method correctly chooses the dimension under different models.

\begin{figure}[b]
  \centering
  \includegraphics[width=7cm]{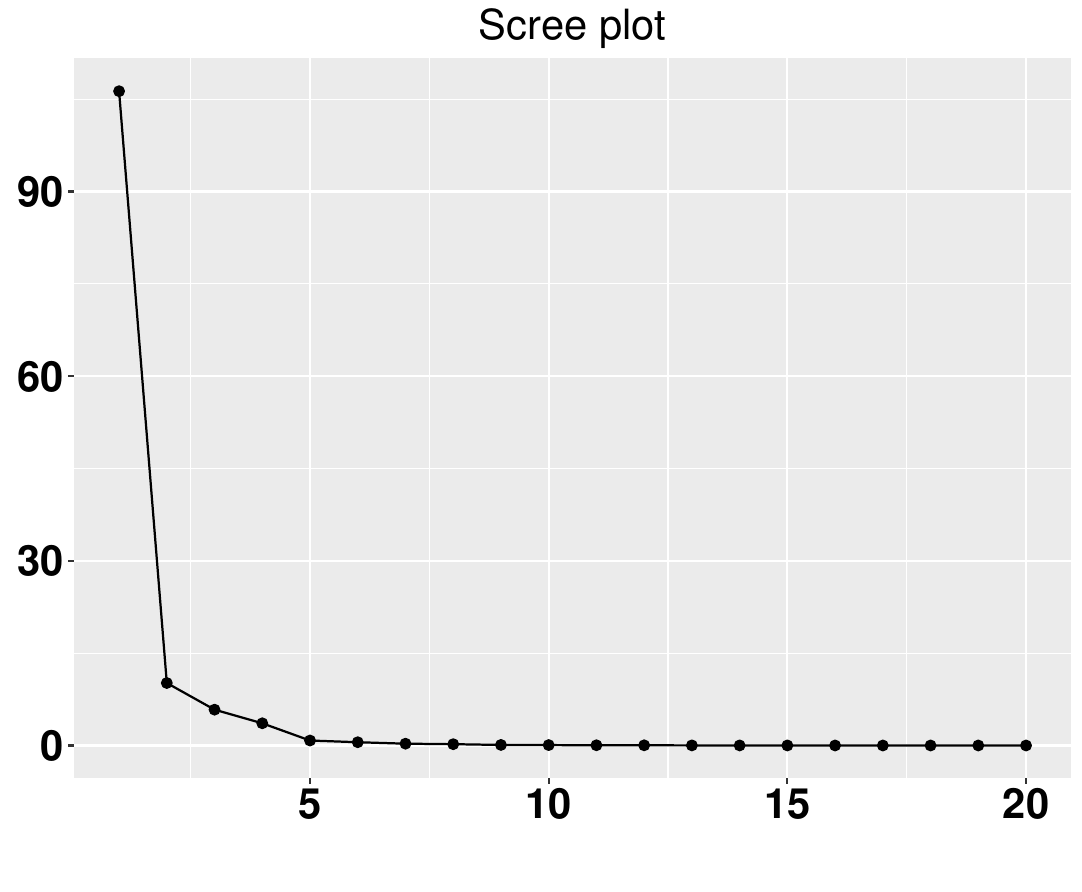}
  \caption{Eigenvalues of near infrared absorbance spectrum data.}\label{scree}
\end{figure}

\subsection{Real data analysis}
We consider the Tecator spectrometric data, available at \url{http://lib.stat.cmu.edu/datasets/tecator} and \texttt{R} package \texttt{fda.usc}. These data are recorded on a Tecator Infratec Food and Feed Analyzer
working in the wavelength range 850 - 1050 nm by the Near Infrared Transmission (NIT) principle.
For each meat sample the data consists of a 100 channel spectrum of
absorbance and the contents of moisture, fat and protein.
The absorbance is $-\log_{10}$ of the transmittance
measured by the spectrometer. The three contents, measured in percent,
are determined by analytic chemistry.
In this example, the task is to predict the fat content $U$ of a
meat sample on the basis of its near infrared absorbance spectrum $X$.
The spectral data $X$ is the functional predictor and fat $U$ is the scalar variable.
In accordance with the literature \cite{ferre2005smoothed},\cite{hsing2009rkhs}, we use the transformed  $Y=\log_{10}(U/(1-U))$ as the response.
\begin{figure}[t]
\centering
\subfigure{
\includegraphics[width=3.5cm]{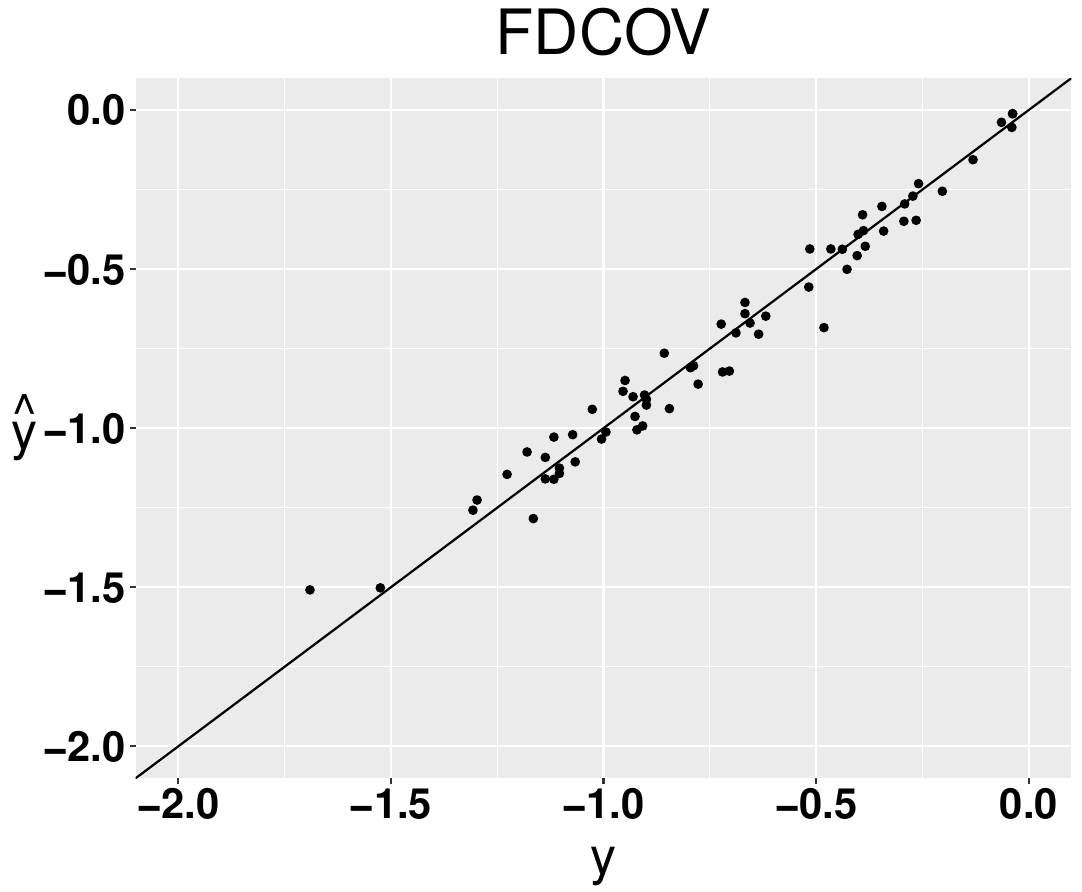}
}
\quad
\subfigure{
\includegraphics[width=3.5cm]{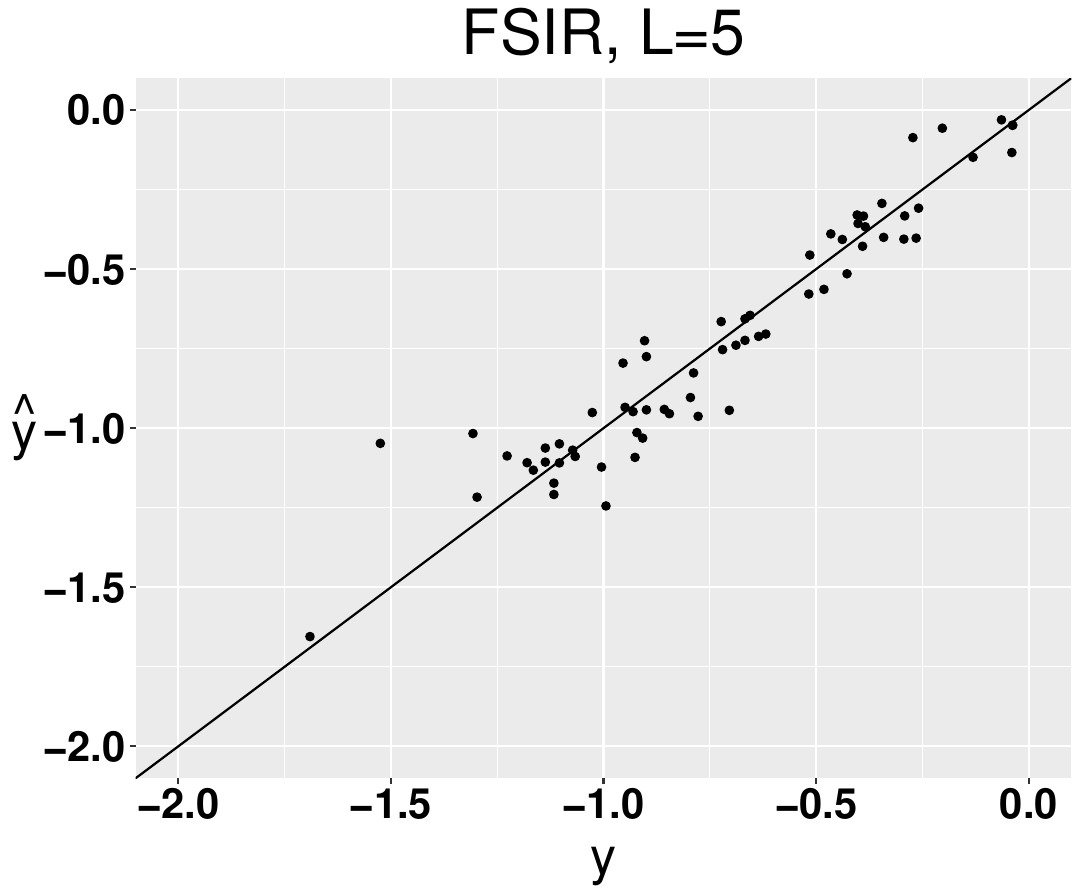}
}
\quad
\subfigure{
\includegraphics[width=3.5cm]{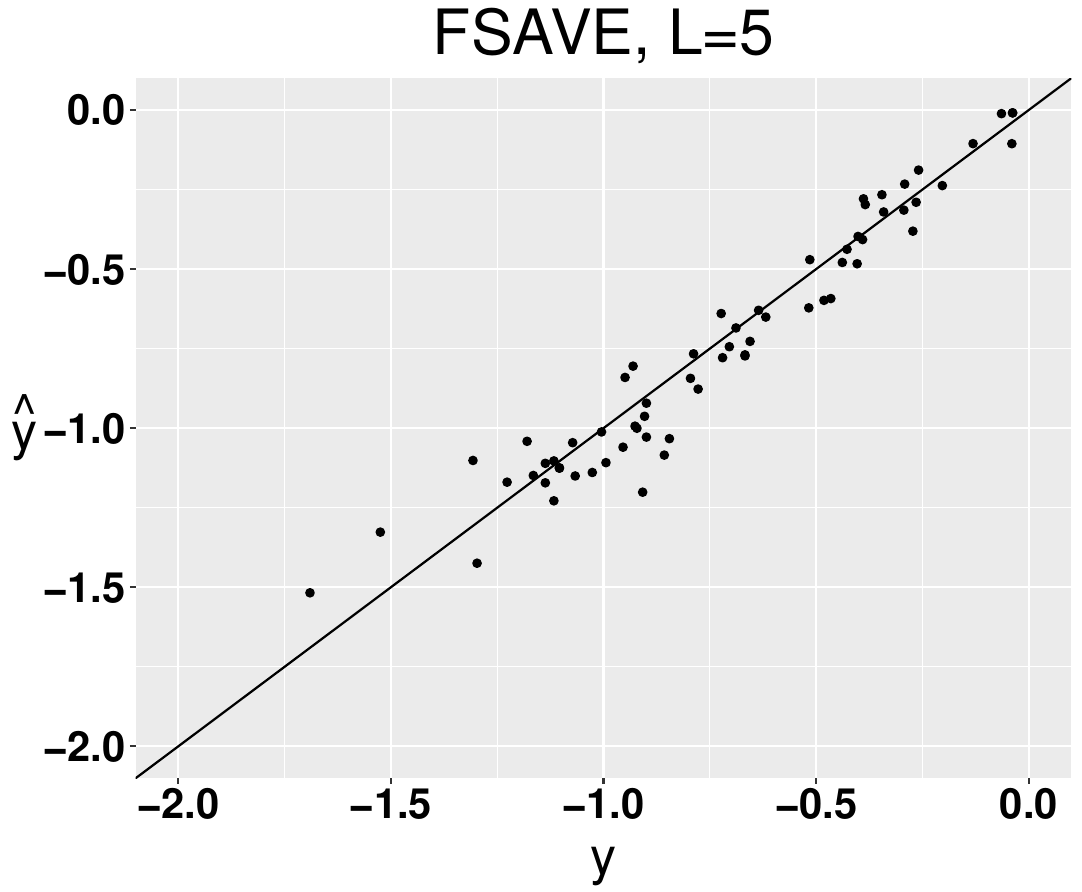}
}
\\
\subfigure{
\includegraphics[width=3.5cm]{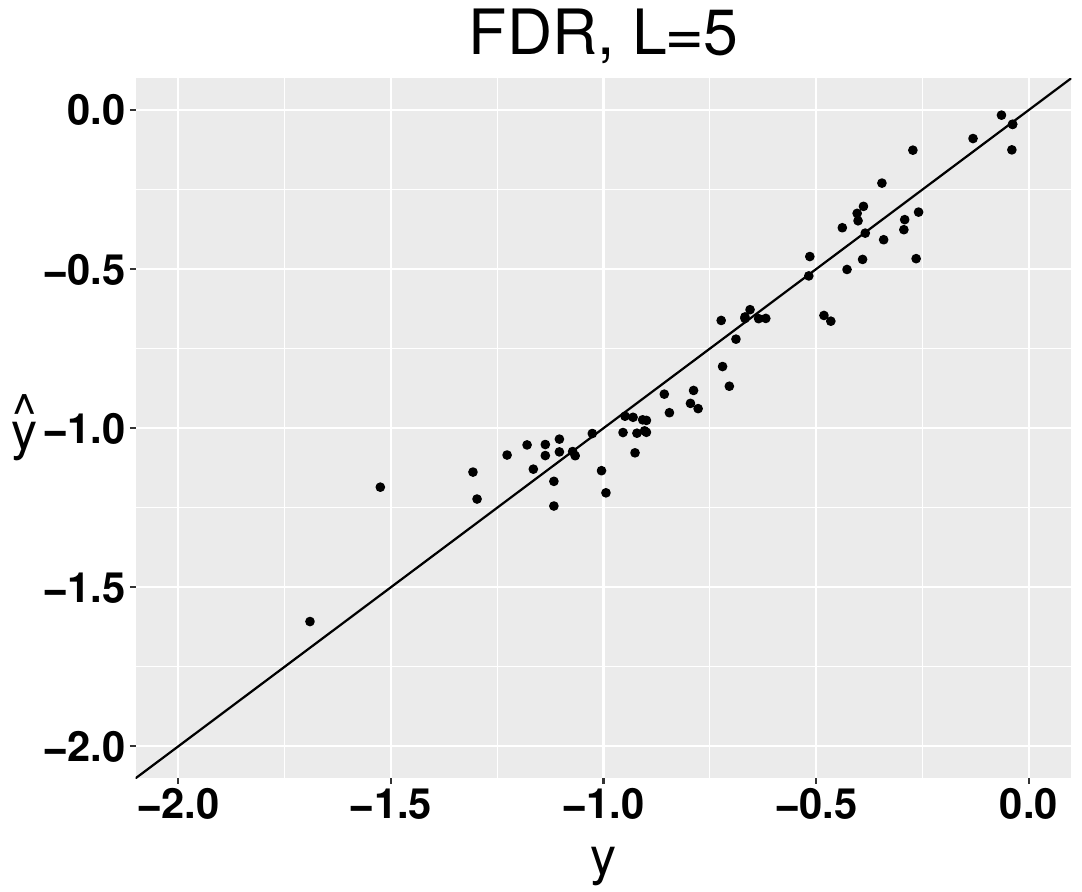}
}
\quad
\subfigure{
\includegraphics[width=3.5cm]{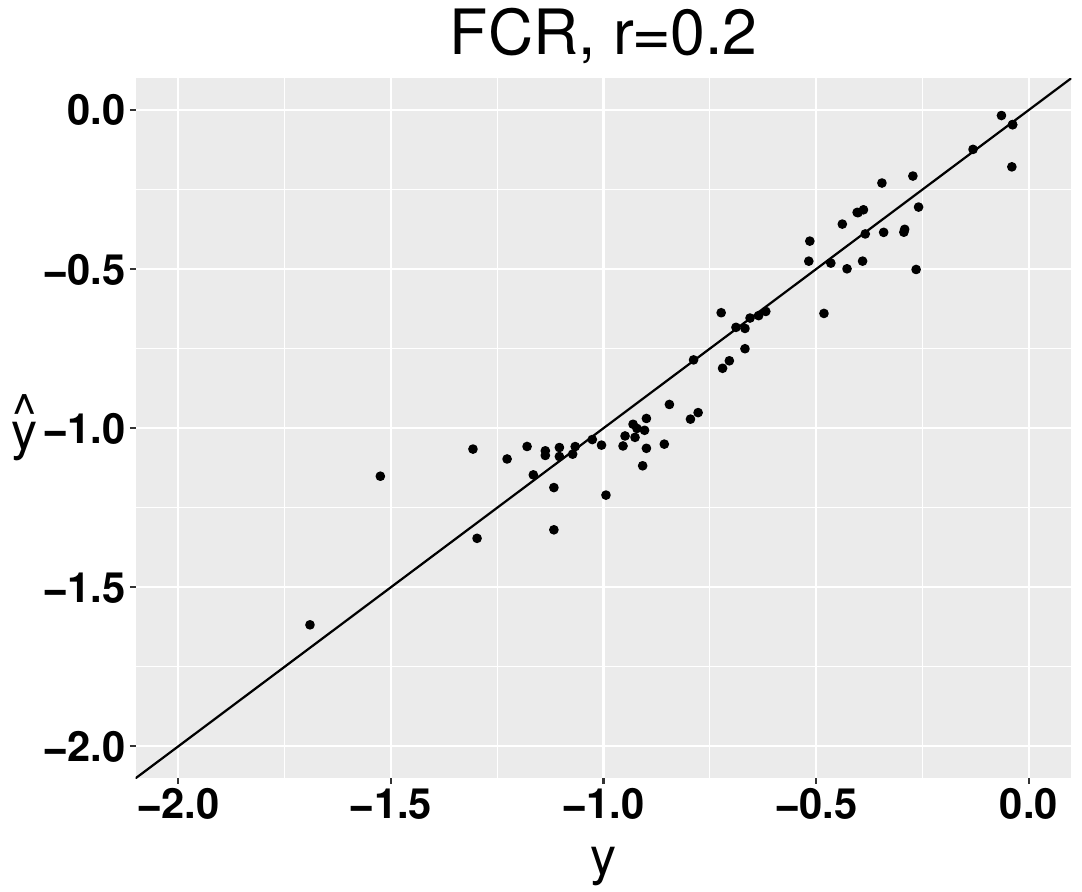}
}
\caption{$\hat{y}$ versus $y$ for different methods.}\label{realdatafig}
\end{figure}

\begin{table}[t]
\centering
\begin{threeparttable}[t]
\caption{Prediction errors  for different estimators for real data. }\label{Realdata}
\begin{tabular}{cccccccc}
\hline
    & $L$& $r$  & FDCOV & FSIR  & FSAVE & FDR  & FCR      \\ \hline
    & 5  & 0.05 & 0.094 & 0.135 & 0.150 & 0.154& 0.147 \\
    & 10 & 0.10 & 0.094 & 0.132 & 0.148 & 0.156& 0.159    \\
    & 15 & 0.15 & 0.094 & 0.127 & 0.153 & 0.159& 0.146   \\
    & 20 & 0.20 & 0.094 & 0.133 & 0.133 & 0.169& 0.137   \\ \hline
\end{tabular}
\end{threeparttable}
\end{table}
The sample size of the data is $n=215$ and we use the first 150 for training and the remaining 65 for testing.
We perform a spectral decomposition of $X$ and draw a scree plot of the eigenvalues in Figure \ref{scree}.
The scree plot shows that first $5$ eigenvectors explain almost the total variation and the first $5$
eigenvalues are not ``too small". Therefore, we select $5$ as the truncation number.
We apply all the methods in the simulation to this real  data and we also consider $L=5,\ 10,\ 15,\ 20$ and $r=0.05,\ 0.10,\ 0.15,\ 0.20$ and use the bootstrap method mentioned in Subsection \ref{subsecK} to select the dimension $K$.
The estimated dimension is $\widehat{K}=4$ in all cases. For every method, we obtain $4$ estimated  projections
$\big(\widehat{\xi}_1,\ \widehat{\xi}_2,\ \widehat{\xi}_3,\ \widehat{\xi}_4\big)=
\big(\langle\widehat{\eta}_1,X\rangle,\ \langle\widehat{\eta}_2,X\rangle,\ \langle\widehat{\eta}_3,X\rangle,\ \langle\widehat{\eta}_4,X\rangle\big)$ to estimate the unknown link function $g$. We use smoothing spline ANOVA method
(\texttt{ssanova} function in the \texttt{R} package \texttt{gss} \citep{gu2013smoothing}).  To measure the predictive performance of different methods, we use the root mean squared prediction error (RMSE) of the test sample which is defined as
\begin{equation}\label{}
  \text{RMSE}= \sqrt{n_{te}^{-1}\sum_{j=1}^{n_{te}}(y_j-\widehat{y}_j)^2},
\end{equation}
where $n_{te}=65$ is the size of test sample, $\widehat{y}_j$ is the predictive value and $y_j$ is the corresponding observed value.
We also use 5-fold cross-validation discussed in Subsection \ref{subsecK} to select parameters $D$ and $K$.
The parameters selected by the cross-validation method are $(D,K)=(5,4)$,  the same as those of the previous method.
The results of the prediction errors for different estimators are reported in Table \ref{Realdata}. From Table \ref{Realdata}, our proposed method outperforms FSIR,
FSAVE, FDR and FCR under all settings.
Figure \ref{realdatafig} is the plot of $\widehat{y}=\widehat{g}\big(\widehat{\xi}_1,\ \widehat{\xi}_2,\ \widehat{\xi}_3,\ \widehat{\xi}_4\big)$ versus $y$ for our proposed  FDCOV method and other competing methods on the test sample.
The figure shows that the predictive responses for the proposed method are really closer to the test sample
than some other methods,
which indicates that our proposed  method retains enough information for regression to predict the response variable.

\section{Concluding remarks}\label{condis}
In this work we propose a method of sufficient dimension reduction for functional data using distance covariance and
establish its statistical consistency. In the estimation procedure, we adopt the commonly used functional PCA to project the infinite-dimensional predictor onto a finite-dimensional subspace. We develop the FDCOV method to estimate $S_{y|x}$, along with procedures for determining the structural dimension.
FDCOV requires very mild conditions on the predictor, unlike the existing methods require the restrictive linear conditional
mean assumption and constant covariance assumption. It also does not involve the inverse of the covariance operator
which is not bounded. In addition, the proposed method does not need to tune the parameters but other methods are needed. For example, the number of slices in FSIR, FSAVE and FDR and the proportion of empirical directions in FCR.

In practice use, other basis such as wavelets and B-spline can also be considered for projection. The theoretical properties such as convergence rate and  asymptotic normality need also be established. We consider the case that the response is a scalar in this article. However, the response can also be a random vector \citep{chen2019sufficient} or a random function \citep{li2017nonlinear} and  the method will be adjusted accordingly in these cases. Nonlinear functional  sufficient dimension reduction methods can also be developed by means of RKHS \citep{li2017nonlinear}.
We leave these to future work.

\section{Appendix}
\emph{Proof of Proposition \ref{prop1}.} Since $span(\bm{\beta})\subseteq span(\bm{\eta})=S_{y|x}$, $K_1\leq K$ we can find a $K\times K_1$ matrix $A$, which satisfies $\bm{\beta}=\bm{\eta}A$. Therefore,
$\mathcal{V}^2(\langle \bm{\beta},X\rangle,Y)=\mathcal{V}^2(A^T\langle \bm{\eta},X\rangle,Y)$. Assume the single value decomposition of $A$ is $U\Lambda V^T$, where $U$ is a $K\times K$ orthogonal matrix, $V$ is a $K_1\times K_1$ orthogonal matrix, and $\Lambda$ is a $K\times K_1$ diagonal matrix. Since $I_{K_1}=\langle\bm{\beta} , \bm{\beta}\rangle_{\Sigma_X}=\langle\bm{\eta}A , \bm{\eta}A\rangle_{\Sigma_X}=A^T\langle\bm{\eta} , \bm{\eta}\rangle_{\Sigma_X}A=A^TA$, we have that all nonzero numbers on the diagonal of $\Sigma$ are $1$. According to the property (ii) in Subsection \ref{DC},
$\mathcal{V}^2(\langle \bm{\beta},X\rangle,Y)=\mathcal{V}^2( V\Lambda^T U^T\langle\bm{\eta},X\rangle,Y)
=\mathcal{V}^2( \Lambda^T U^T\langle\bm{\eta},X\rangle,Y)$.

Denote $U^T\langle\bm{\eta},X\rangle=(Z_1,\ldots,Z_K)^T$. Since all nonzero numbers on the diagonal of $\Lambda$ are $1$, we have $\Lambda^TU^T\langle\bm{\eta},X\rangle=(Z_1,\ldots,Z_{K_1})^T$. Clearly,  $\Lambda^TU^T\langle\bm{\eta},X\rangle$ is a  vector composed of the first $K_1$ components of $U^T\langle\bm{\eta},X\rangle$.
From this observation and Lemma A.1 in \cite{sheng2016sufficient}, we have $\mathcal{V}^2( \Lambda^T U^T\langle\bm{\eta},X\rangle,Y)\leq\mathcal{V}^2( U^T\langle\bm{\eta},X\rangle,Y)$ and the equality holds if and only if $K_1=K$. By property (ii) in Subsection \ref{DC},
$\mathcal{V}^2( U^T\langle\bm{\eta},X\rangle,Y)\leq\mathcal{V}^2( \langle\bm{\eta},X\rangle,Y)$.
Thus, we obtain $\mathcal{V}^2( \langle\bm{\beta},X\rangle,Y)\leq\mathcal{V}^2( \langle\bm{\eta},X\rangle,Y)$, and the equality holds if and only if $span(\bm{\beta})=span(\bm{\eta})$. $\hfill\blacksquare$

\vskip0.5cm

\emph{Proof of Proposition \ref{prop2}.} For the $\bm{\beta}$ and $\bm{\eta}$ defined in Proposition \ref{prop2}, we can find a rotation matrix $R$ such that $\bm{\beta}R=(\bm{\eta}_a, \bm{\eta}_b)$ and
$span(\bm{\eta}_a)\subseteq span(\bm{\eta})$, $span(\bm{\eta}_b)\subseteq span(\bm{\eta})^{\perp}$ where $span(\bm{\eta})^{\perp}$ is the  orthogonal complement space of $span(\bm{\eta})$ with respect to the inner product $\langle\cdot ,\cdot \rangle_{\Sigma_X}$.

The definition of $\bm{\eta}$ indicates $Y\perp X|\langle\bm{\eta},X\rangle$, thus
$Y\perp \langle \bm{\eta}_b,X\rangle|\langle\bm{\eta},X\rangle$. By  Condition \ref{assump1}, we have
$\langle\bm{\eta}_b, X\rangle\perp\langle\bm{\eta}, X\rangle$.
Therefore $\Big(\begin{array}{c}
            Y \\
            \langle\bm{\eta}, X\rangle
          \end{array}\Big) \perp \langle\bm{\eta}_b, X\rangle$, and we can get
          $\Big(\begin{array}{c}
            Y \\
            \langle\bm{\eta}, X\rangle
          \end{array}\Big) \perp \langle\bm{\eta}_b, X\rangle$ by  Proposition 4.3 in \cite{cook2009regression}.
          Let $U_1=\Big(\begin{array}{c}
            \langle\bm{\eta}_a, X\rangle\\
            \bm{0}
          \end{array}\Big) $, $V_1=Y$, $U_1=\Big(\begin{array}{c}
            \bm{0}\\
            \langle\bm{\eta}_b, X\rangle
          \end{array}\Big) $ and $V_2=0$, then $(U_1,V_1)\perp (U_2,V_2)$.
          By property (iii) in  Subsection \ref{DC}, $\mathcal{V}^2(U_1+U_2,V_1+V_2)<\mathcal{V}^2(U_1+V_1)+\mathcal{V}^2(U_2+V_2)$,
          this means $\mathcal{V}^2(R^T\langle\bm{\beta}, X\rangle,Y)=\mathcal{V}^2(\langle\bm{\beta}, X\rangle,Y)
          <\mathcal{V}^2(\langle\bm{\eta}_b, X\rangle,Y)\leq \mathcal{V}^2(\langle\bm{\eta}, X\rangle,Y)$. $\hfill\blacksquare$

\vskip0.5cm
\emph{Proof of Theorem \ref{thm1}.} Suppose $\bm{\eta}_{n,D_n}$ is not a consistent estimator of $S_{y|x}$,
there exists a subsequence $\bm{\eta}_{n^*,D_{n^*}}$ of $\bm{\eta}_{n,D_n}$ such that
$\bm{\eta}_{n^*,D_{n^*}} \stackrel{P}{\longrightarrow} \bm{\eta}^*$, where $\langle\bm{\eta}^* ,\ \bm{\eta}^*\rangle_{\widehat{\Sigma}_X}=I_{K}$ but $span(\bm{\eta}^*)\neq span(\bm{\eta})$. By Lemma A in \cite{sheng2016sufficient}, we have
$$\mathcal{V}_n^2(\langle \bm{\eta}_{n^*,D_{n^*}},\bm{X}\rangle,\bm{Y})-
\mathcal{V}_n^2(\langle\bm{\eta}^*,\bm{X}\rangle,\bm{Y})\stackrel{P}{\longrightarrow}0.$$ According to Theorem 2 in \cite{szekely2007measuring}, $\mathcal{V}_n^2(\langle\bm{\eta}^*,\bm{X}\rangle,\bm{Y})\stackrel{a.s.}{\longrightarrow}
\mathcal{V}^2(\langle\bm{\eta}^*,\bm{X}\rangle,\bm{Y})$, therefore $\mathcal{V}_n^2(\langle \bm{\eta}_{n^*,D_{n^*}},\bm{X}\rangle,\bm{Y})\stackrel{P}{\longrightarrow}\mathcal{V}^2(\langle\bm{\eta}^*,\bm{X}\rangle,\bm{Y})$.

Besides, since   $\bm{\eta}_{n,D_n}=\arg\max_{\langle\bm{\beta}^{D_n} ,\ \bm{\beta}^{D_n}\rangle_{\widehat{\Sigma}_X}=I_{K}} \mathcal{V}_n^2(\langle \bm{\beta}^{D_n},\bm{X}\rangle,\bm{Y})$, we have
$$\mathcal{V}_n^2(\langle \bm{\eta}_{n,D_n},\bm{X}\rangle,\bm{Y})\geq \mathcal{V}_n^2(\langle \bm{\eta}^{D_n},\bm{X}\rangle,\bm{Y}),$$ where  $\bm{\eta}^{D_n}$ is the representation of the function $\bm{\eta}$ in the
$D_n$-truncated basis.
Let $n\rightarrow \infty$, we get  $\mathcal{V}^2(\langle \bm{\eta}^*,\bm{X}\rangle,\bm{Y})\geq \mathcal{V}^2(\langle \bm{\eta},\bm{X}\rangle,\bm{Y})$, which contradicts to the definition of $\bm{\eta}$.
Then we can conclude that $\bm{\eta}_{n,D_n}$ is a consistent estimator of $\bm{\eta}$.
$\hfill\blacksquare$

\bibliographystyle{tfnlm}
\bibliography{fSDRviaDCOV}
\end{document}